\documentclass{amsart}

\usepackage{amsmath,amsfonts,amsthm,amssymb}
\usepackage{graphicx}

\newcommand{\mapdown}[1]%
{\Big\downarrow\rlap{$\vcenter{\hbox{$\scriptstyle#1$}}$}}

\newcommand{\A}{\mathcal{A}}
\newcommand{\arr}{\longrightarrow}
\newcommand{\spec}[1]{\mathop{Sp}\left(#1\right)}
\newcommand{\symm}[1][\alb]{\mathop{\mathrm{Symm}}\left(#1\right)}
\newcommand{\supp}{\mathop{\mathrm{supp}}}
\newcommand{\C}{\mathbb{C}}
\newcommand{\R}{\mathbb{R}}
\newcommand{\st}{\mathop{\mathrm{St}}}
\newcommand{\wt}{\widetilde}
\newcommand{\alb}{\mathsf{X}}
\newcommand{\xo}{\alb^\omega}
\newcommand{\xmo}{\alb^{-\omega}}
\newcommand{\xs}{\alb^*}
\newcommand{\xz}{\alb^{\mathbb{Z}}}
\newcommand{\maxalg}{\mathcal{A}_{\textit{max}}}
\newcommand{\minalg}{\mathcal{A}_{\textit{min}}}
\newcommand{\nuke}{\mathcal{N}}
\newcommand{\cpg}{O_G}
\newcommand{\solen}[1][G]{\mathcal{S}_{#1}}
\newcommand{\si}{\widehat{\mathsf{s}}}
\newcommand{\sis}{\mathsf{s}}
\newcommand{\lims}[1][G]{\mathcal{J}_{#1}}
\newcommand{\Z}{\mathbb{Z}}
\newcommand{\img}[1]{\mathop{\mathrm{IMG}}\left(#1\right)}

\newcommand{\M}{\mathcal{M}}
\newcommand{\gi}{\M_{d^\infty}(G)}

\newtheorem{theorem}{Theorem}[section]
\newtheorem{proposition}[theorem]{Proposition}
\newtheorem{corollary}[theorem]{Corollary}
\newtheorem{lemma}[theorem]{Lemma}

\theoremstyle{definition}

\newtheorem{definition}{Definition}[section]
\newtheorem*{remark}{Remark}
\newtheorem{example}{Example}
\newtheorem{examples}[example]{Examples}
\newtheorem{question}{Problem}

\numberwithin{equation}{section}

\title{Self-similar groups, operator algebras and Schur complements}
\author{Rostislav Grigorchuk}
\address{Texas A\&M University, College Station, USA}
\email{grigorch@math.tamu.edu}
\thanks{Supported by NSF grants DMS-0456185 and DMS-0600975}

\author{Volodymyr Nekrashevych}
\email{nekrash@math.tamu.edu}
\thanks{Supported by NSF grant DMS-0605019}
\dedicatory{Dedicated to the 70th birthday of D.V.~Anosov}

\begin{document}
\maketitle

\begin{abstract}
In the first part of the article we introduce $C^*$-algebras
associated to self-similar groups and study their properties and
relations to known algebras. The algebras are constructed as
sub-algebras of the Cuntz-Pimsner algebra (and its homomorphic
images) associated with the self-similarity of the group. We study
such properties as nuclearity, simplicity and Morita equivalence
with algebras related to solenoids.

The second part deals with the Schur complement transformations of
elements of self-similar algebras. We study properties of such
transformations and apply them to the spectral problem for Markov
type elements in self-similar $C^*-$algebras. This is related to
the spectral problem of the discrete Laplace operator on groups
and graphs. Application of the Schur complement method in many
situations reduces the spectral problem to study of invariant sets
(very often of the type of a ``strange attractor'') of a
multidimensional rational transformation. A number of illustrating
examples is provided. Finally we observe a relation between the
Schur complement transformations and Bartholdi-Kaimanovich-Virag
transformations of random walks on self-similar groups.
\end{abstract}

\section{Introduction}

Self-similar groups is a class of groups which attracts more and
more attention of researchers from different areas of mathematics,
and first of all from group theory.

Self-similar groups (whose study was initiated by the first named
author, N.~Gupta, S.~Sidki, A.~Brunner, P.M.~Neumann, and others)
posses many nice and unusual properties which allow to solve
difficult problems of group theory and related areas, even
including problems in Riemannian geometry and holomorphic
dynamics. They are closely related to groups generated by finite
automata (studied by J.~Ho{\v r}ej{\v s}, V.M.~Glushkov, S.
V.~Aleshin, V.I.~Sushchanski and others) and many of them belong
to another interesting class of groups---the class of branch
groups. We recommend the following sources for the basic
definitions and properties
\cite{nek:book,handbook:branch,grigorchuk:branch,sidki_monogr}.

One of the main features that makes the class of self-similar
groups important is that it makes possible treatment of the renorm
group in a noncommutative setting. This passage from cyclic to a
non-commutative renormalization can be compared with the passage
from classical to non-commutative geometry, i.e., passage from
commutative $C^*$-algebras of continuous functions to
non-commutative $C^*$-algebras, see~\cite{connes:noncomg}.

Recent research shows that operator algebras also play an
important role in the theory of self-similar groups and show
interesting connections with other areas (for instance with
hyperbolic dynamics). The first appearance of $C^*-$algebras
related to self-similar groups was in~\cite{bgr:spec} and was
related to the problem of computation of the spectra of Markov
type operators on the Schreier graphs related to self-similar
groups. In~\cite{nek:bim} the methods of~\cite{bgr:spec} were
interpreted in terms of the Cuntz-Pimsner algebras of the Hilbert
bimodules associated with self-similar groups. It was proved in~\cite{nek:bim}
that there exists a smallest self-similar norm (i.e., norm which
agrees with the bimodule) and that the Cuntz-Pimsner
algebra constructed using the completion of the group algebra by
the smallest norm is simple and purely infinite.

This article is a survey of results and ideas on the interplay of
self-similar groups, renormalization, self-similar operator
algebras, spectra of Markov operators and random walks.

We continue the study of self-similar completions of the group
algebra started in~\cite{bgr:spec} and~\cite{nek:bim}. We prove
existence of the largest self-similar norm on the group algebra
and define the maximal Cuntz-Pimsner algebra and maximal
self-similar completion of the group algebra
(Propositions~\ref{pr:cpimsner} and~\ref{pr:algmaximal}).

Our starting point is the observation that self-similarities of a
Hilbert space $H$ (i.e., isomorphisms between $H$ and its $d-$th,
$d\ge 2$ power $H^d$) are in a natural bijection with the
representations of the Cuntz algebra $O_d$. Then we consider
self-similar unitary representations of a self-similar group and
following~\cite{nek:bim} define the associated universal
Cuntz-Pimsner algebra $\cpg$. The algebra $\cpg$ is universal in
the sense that any unitary self-similar representation of the
group can be extended to a representation of $\cpg$. Homomorphic
images of $\cpg$ can be constructed using different self-similar
representations of $G$. Among important self-similar
representations we study the natural unitary representation of $G$
on $L^2(\xo, \nu)$, where $\xo$ is the boundary of the rooted tree
on which $G$ acts and $\nu$ is the uniform Bernoulli measure on
it. Another important class of self-similar representations are
permutational representations of $G$ on countable $G$-invariant
subsets of $\xo$.

The sub-$C^*$-algebra of $\cpg$ generated by $G$ is denoted
$\maxalg$. For every homomorphic image of the Cuntz-Pimsner
algebra $\cpg$ we get the respective image of $\maxalg$. The
algebra $\minalg$, for instance, is defined as the image of
$\maxalg$ under a permutational representation of $G$ on a generic
self-similar $G$-invariant countable set of $\xo$ (which can be
extended to a permutational representation of $\cpg$ in a natural
way). The algebra $\minalg$ was studied in~\cite{bgr:spec}
and~\cite{nek:bim}. Similarly, the algebra $\A_{mes}$ generated by
the natural representation of $G$ (and of $\cpg$) on $L^2(\xo,
\nu)$, is considered. This algebra is particularly convenient for
spectral computations (see Subsection~\ref{ss:scssalg}).

It is convenient in many cases to pass to a bigger sub-algebra of
$\cpg$ than $\maxalg$. It is the algebra generated by $G$ and the
union of the matrix algebras $M_{d\times d}(\C)$ naturally
constructed inside $O_d\subset\cpg$. It is proved
in~\cite{nek:cpalg} that this algebra is denoted
$\M_{d^\infty}(G)$ and it is the universal algebra of the groupoid
of germs of the action of $G$ on the boundary $\xo$ of the rooted
tree. This makes it possible to apply the well-developed theory of
$C^*$-algebras associated to groupoids to the study of $\cpg$ and
$\maxalg$.

For every homomorphic image of $\cpg$ (i.e., for every
self-similar representation of $G$) we get the corresponding image
of $\M_{d^\infty}(G)$.

In Section~\ref{s:contracting} we show a relation of self-similar
groups and algebras to hyperbolic dynamics. If a self-similar
group $G$ is \emph{contracting}, then there is a natural Smale
space associated to it (the \emph{limit solenoid}). It is a
dynamical system $(\solen, \si)$ with hyperbolic behavior: the
space $\solen$ has a local structure of a direct product such that
the homeomorphism $\si$ is contracting on one factor and expanding
on the other. We prove in Theorem~\ref{th:moritaequiv} that the
algebra $\M_{d^\infty}(G)$ is Morita equivalent in this case to
the convolution algebra of the unstable equivalence relation on
the limit solenoid. Such convolution algebras were studied by
J.~Kaminker, I.~Putnam and J.~Spielberg
in~\cite{putnam,put_kam_sp,putnam:structure}.

The algebras $\maxalg$, $\M_{d^\infty}(G)$ and their images under
self-similar representations of $G$ have a nice self-similarity
structure described by matrix recursions, which encode the
structure of the Moor diagrams of the Mealy type automata defining
the underlying group. This self-similarity is used in the second
part of the article in the study of spectra of elements of the
involved algebra.

We recall at first the classical Schur complement transformation,
which is used in linear algebra for solving systems of linear
equations, in statistics for finding conditional variance of
multivariate Gaussian random variables and in Bruhat normal form
(see~\cite{cottle:schursurvey,cohn}). We show in our article that
Schur complement is also useful in the study of spectral problems
in self-similar algebras and that it can be nicely expressed in
terms of the Cuntz algebra. We establish some simple properties of
the Schur complement transformations and introduce a semigroup of
such transformations.

The method that we use to treat the spectral problem could be a
first step in generalization of the method of Malozemov and
Teplyaev that they developed for the study of spectra of
self-similar graphs related to classical fractals (like the
Sierpinskii gasket)~\cite{teplyaev:sierpinski,malozemov_teplyaev}.
In fact they also use (in an implicit form) the Schur complement.
Their technique is developed for the case when only one complex
parameter is involved. Our technique (which is a development of
the technique used in~\cite{bgr:spec,gr_zu:lamp,grisunik:hanoi}
involves several parameters and therefore necessarily lead to
multidimensional rational mappings and their dynamics.

The Schur complements in our situations are renormalization
transformations for the spectral problem and related problems.
Considered together they generate a noncommutative semigroup which
can be called the ``Schur renorm group'' (observe that in
classical situation the renorm group very often is a cyclic
semigroup).

We illustrate our method by several examples the most
sophisticated among which is the example related to the
3-generated torsion 2-group of intermedate growth constructed
in~\cite{grigorchuk:80_en}.

The transformations that arise in this case are
\[\wt S_1:\left(\begin{array}{c} x\\ y\\ z\\ u\\
v\end{array} \right)\mapsto\left(\begin{array}{c} z+y\\
\frac{x^2(2yzv-u(y^2+z^2-u^2+v^2))}{(y+z+u+v)(y+z-u-v)(y-z+u-v)(-y+z+u-v)}\\
\frac{x^2(2zuv-y(-y^2+z^2+u^2+v^2))}{(y+z+u+v)(y+z-u-v)(y-z+u-v)(-y+z+u-v)}\\
\frac{x^2(2yuv-z(y^2-z^2+u^2+v^2))}{(y+z+u+v)(y+z-u-v)(y-z+u-v)(-y+z+u-v)}\\
u+v+\frac{x^2(2yzu-v(y^2+z^2+u^2-v^2))}{(y+z+u+v)(y+z-u-v)(y-z+u-v)(-y+z+u-v)}\end{array}\right)\]
and
\[
\wt S_2:\left(\begin{array}{c} x\\ y\\ z\\ u\\
v\end{array}\right)\mapsto\left(\begin{array}{c}\frac{x^2(y+z)}{(u+v+y+z)(u+v-y-z)}\\
u\\ y\\ z\\
v-\frac{x^2(u+v)}{(u+v+y+z)(u+v-y-z)}\end{array}\right)
\]

The dynamical properties of these transformations are not well
understood and this is one of intriguing problems.

The examples that appear demonstrate few cases (taken from
\cite{bgr:spec,grisunik:hanoi}) leading to easily treatable
transformations (when they ``integrable'' in the sense that there
is a nontrivial semi-conjugacy to a one-dimensional map), a couple
of examples (taken from
\cite{zukgrigorchuk:3st,grisavchuksunic:img}) when there is no
information about topological nature of the invariant subsets
(which according to computer experiments look like ``strange
attractors''), and a couple of examples when our method does not
work (but the corresponding spectral problem is important because
of its relation to the problem of finding new constructions of
expanders).

In principle, the area of application of our method is much
broader but then it requires use of Schur complements in
infinite-dimensional spaces of matrices. At the moment we do not
have examples of successful applications of the method in the
infinite-dimensional case.

One of important properties we are after in the study of
self-similar groups and related objects is amenability.
Self-similar groups provide a number of examples of amenable but
not elementary amenable groups
\cite{grigorchuk:growth_en,grigorchuk:notEG}. The fundamental idea
how to treat amenability of self-similar groups belongs to
Bartholdi and Virag~\cite{barthvirag}. Roughly speaking it
converts self-similarity of a group into self-similarity of a
random walk on it. This idea was further developed by
V.Kaimanovich (in what he calls the ``M\"unchhausen trick'') using
the notion of entropy of a random
walk~\cite{kaimanovich:munchhausen}. Kaimanovich introduced
transformations of the random walk under the renorm
transformations of the group and successfully used them to show
that some self-similar groups are amenable. Our observation is
that these transformations again can be interpreted as the Schur
complement transformations of the measures (or corresponding
elements of the group algebra) determining the random walk. More
precisely, we show that the transformation considered
in~\cite{kaimanovich:munchhausen} is the Schur complement
conjugated by the map $A\mapsto A+I$, where $I$ is the identity
matrix.

We believe that the introduced ideas of self-similar algebras and
Schur complement transformations on them will be useful for the
study of different aspects of the theory of self-similar groups
and its applications.

\medskip\noindent\textbf{Acknowledgements.} The authors are
grateful to Peter Kuchment for pointing out that the
transformation we are considering is called the Schur complement
and to Yaroslav Vorobets for useful remarks.

\section{Self-similar groups}
\subsection{Definition}
Let $\alb$ be a finite alphabet and let $\xs$ denote the set of
finite words over the alphabet $\alb$. In other terms, $\xs$ is
the free monoid generated by $\alb$. We consider $\xs$ to be the
set of vertices of a rooted regular tree with the root coinciding
with the empty word $\varnothing$ and in which a word $v$ is
connected to every word of the form $vx$ for $x\in\alb$.

\begin{figure}
 \includegraphics{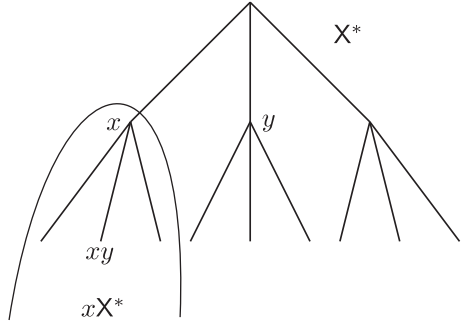}\\
 \caption{Rooted tree}\label{fig:tree}
\end{figure}

A rooted tree is a standard model of a self-similar structure. If
we remove the empty word from it, then the tree $\xs$ will split
into $|\alb|$ subtrees $x\xs$, where $x\in\alb$. Each of the
subtrees $x\xs$ is isomorphic to the whole tree $\xs$ under the
isomorphism $xv\mapsto v$, see Figure~\ref{fig:tree}.

If we consider the \emph{boundary} of the tree $\xs$, then the
self-similarity is even more evident. The boundary of a rooted
tree is the set of all infinite paths starting in the root. In our
case the boundary of $\xs$ is naturally identified with the set
$\xo$ of infinite words $x_1x_2\ldots$. The set $\xo$ is a
disjoint union of the \emph{cylindrical sets}
$x\xo=\{xx_2x_3\ldots\;:\;x_i\in\alb\}$, and again, the shift
$xw\mapsto w$ is a bijection (and a homeomorphism, if we endow
$\xo$ with the natural topology of a direct product of discrete
sets $\alb$).

A group acting on the rooted tree $\xs$ is called self-similar, if
the action agrees with the described self-similarity structure on
the tree $\xs$. Namely, we adopt the following definition.

\begin{definition}
A \emph{self-similar group} $(G, \alb)$ is a group $G$ acting
faithfully on the rooted tree $\xs$ such that for every $g\in G$
and every $x\in\alb$ there exist $h\in G$ and $y\in\alb$ such that
\[g(xw)=yh(w)\]
for all $w\in\xs$.
\end{definition}

It follows from the definition that if $(G, \alb)$ is a
self-similar group, then for every $v\in\xs$ there exists $h\in G$
uniquely defined by the condition
\[g(vw)=g(v)h(w)\]
for all $w\in\xs$. The element $h$ is called \emph{restriction} of
$g$ in $v$ and is denoted $h=g|_v$.

We have the following obvious properties of restriction
\[g|_{v_1v_2}=\left(g|_{v_1}\right)|_{v_2},\qquad
(gh)|_v=g|_{h(v)}h|_v.\]

It is convenient to identify a letter $x$ with the \emph{creation
operator} on $\xs$ (or on $\xo$) given by appending the letter $x$
to the beginning of the word:
\[x\cdot v=xv.\]

Similarly, we can identify every word $u\in\xs$ with the creation
operator
\[u\cdot v=uv.\]

In this case the condition that $g(vw)=uh(w)$ for all $w\in\xs$
can be written as equality of compositions of transformations of
$\xs$:
\[g\cdot v=u\cdot h.\]

We have the following straightforward corollary of the
definitions.

\begin{proposition}
The set $\alb\cdot G$ of transformations of $\xs$ of the form
$x\cdot g:w\mapsto xg(w)$ is closed under pre- and
post-compositions with action of the elements of $G$. The obtained
left and right actions of $G$ on $\alb\cdot G$ commute and are
defined by
\[h\cdot (x\cdot g)=h(x)\cdot (h|_xg),\qquad (x\cdot g)\cdot
h=x\cdot (gh),\] where dot denotes composition of transformations.
\end{proposition}

We call the set $\alb\cdot G$ the \emph{(permutational)
$G$-bimodule associated with the self-similar group $(G, \alb)$}.

\begin{definition}
\label{def:replicating} A self-similar group $(G, \alb)$ is called
\emph{self-replicating} (\emph{recurrent} in~\cite{nek:book}) if
it is transitive on the first level $\alb^1$ of the tree $\xs$ and
for any (and thus for every) $x\in\alb$ the map $g\mapsto g|_x$
from the stabilizer of $x$ in $G$ to $G$ is onto. Equivalently,
the self-similar action is self-replicating if the left action of
$G$ on the bimodule $\alb\cdot G$ is transitive.
\end{definition}

It is not hard to prove by induction that if an action is
self-replicating, then it is transitive on every level $\alb^n$ of
the tree $\xs$ (is \emph{level-transitive}).

Every self-similar action naturally induces an action on the
boundary $\xo$ of the tree $\xs$. This is an action by
measure-preserving homeomorphisms.

\begin{example} \textbf{Grigorchuk
group.}\label{sss:grigorchuk} Consider the binary alphabet
$\alb=\{0, 1\}$, the tree $\xs$ and let $\mathfrak{G}$ be the
group generated by the automorphisms $a, b, c, d$ of $\xs$ defined
recursively by the relations
\begin{alignat*}{2}
a(0w) &= 1w, &\quad a(1w) &=0w\\
b(0w) &= 0a(w), &\quad b(1w) &= 1c(w)\\
c(0w) &= 0a(w), &\quad c(1w) &= 1d(w)\\
d(0w) &= 0w, &\quad d(1w) &= 1b(w).
\end{alignat*}

\begin{figure}
 \includegraphics{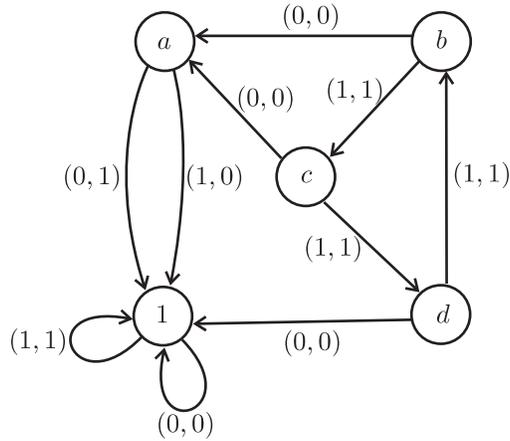}\\
 \caption{Automaton generating the Grigorchuk group}\label{fig:grig}
\end{figure}

This group was defined for the first time
in~\cite{grigorchuk:80_en}. This group is a particularly easy
example of a Burnside group (an infinite finitely generated
torsion group) and it is the first example of a group of
intermediate growth, which answers a question of J.~Milnor.
Figure~\ref{fig:grig} shows the Moore diagram of the automaton
generating the Grigorchuk group (see a definition of Moore
diagrams below).
\end{example}

\begin{example}\textbf{Free group.} \label{sss:free} A convenient way
to define self-similar groups are Moore diagrams of automata,
generating them. Consider, for instance, the Moore diagram shown
on Figure~\ref{fig:alfree}. The vertices of the diagram correspond
to states of the automaton (to generators of the group), the
arrows describe transitions between the states and the labels
correspond to the output of the automaton. If there is an arrow
from a state $g$ to a state $h$ labeled by $(x, y)$, then this
means that $g(xw)=yh(w)$ for all $w\in\xs$. Thus, the automaton
shown on Figure~\ref{fig:alfree} describes the group generated by
the elements $a, b, c$ such that
\begin{alignat*}{2}
a(0w) &= 0b(w), &\quad a(1w)=1b(w),\\
b(0w) &= 1a(w), &\quad b(1w)=0c(w),\\
c(0w) &= 1c(w), &\quad c(1w)=0a(w),
\end{alignat*}
for all $w\in\xs$.

\begin{figure}
 \includegraphics{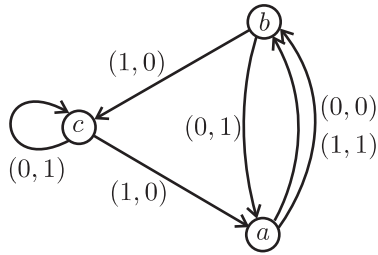}\\
 \caption{Automaton generating a free group}\label{fig:alfree}
\end{figure}

S.~Sidki conjectured in~\cite{sid:cycl} that the three states of
the automaton generate a free group of rank 3. This claim was
later proved by Y.~Vorobets and M.~Vorobets
in~\cite{vorobets:alfree}.
\end{example}

\begin{example}\textbf{Basilica.} \label{sss:basilica} The group
generated by the automaton shown on Figure~\ref{fig:basilica} is
called the \emph{Basilica group}, since it is the iterated
monodromy group of the polynomial $z^2-1$ (see below for the
definition).

\begin{figure}
 \includegraphics{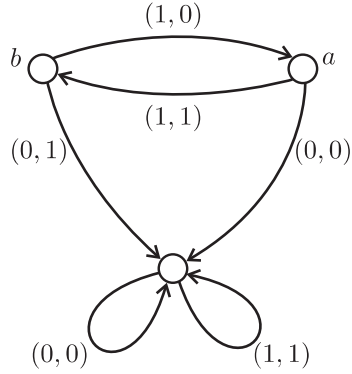}\\
 \caption{Automaton generating $\img{z^2-1}$}\label{fig:basilica}
\end{figure}

It is the first example of an amenable group, which can not be
constructed from the groups of sub-exponential growth using the
operations preserving amenability (taking quotients, extensions,
direct limits and passing to a subgroup).
\end{example}

\subsection{Iterated monodromy groups} Let $\M$ be a path
connected and locally path connected topological space and let
$\M_1$ be its open subset. A \emph{$d$-fold partial self-covering}
is a covering map $f:\M_1\arr\M$, i.e., a continuous map such that
for every $x\in\M$ there exists a neighborhood $U\ni x$ whose
total preimage $f^{-1}(U)$ is a disjoint union of $d$ open subsets
which are mapped homeomorphically onto $U$ by $f$.

For every $n\ge 1$ the iteration $f^n:\M_n\arr\M$ is a $d^n$-fold
partial self-covering (in general with a smaller domain $\M_n$).

Choose a basepoint $t\in\M$. Then the disjoint union
$T=\bigsqcup_{n\ge 0}f^{-n}(t)$ has a natural structure of a
$d$-regular tree with the root $t\in\{t\}=f^{-0}(t)$ where a
vertex $z\in f^{-n}(t)$ is connected to the vertex $f(z)\in
f^{-(n-1)}(t)$.

The fundamental group $\pi_1(\M, t)$ acts naturally on every level
$f^{-n}(t)$ of the tree $T$. The image of a vertex $z$ under the
action of a loop $\gamma\in\pi_1(\M, t)$ is the end of the unique
preimage of $\gamma$ under the covering $f^n$, which starts at
$z$. It is easy to check that these actions define an action of
$\pi_1(\M, t)$ by automorphisms of the rooted tree $T$. This
action is called the \emph{iterated monodromy action}.

\begin{definition}
The \emph{iterated monodromy group} $\img{f}$ of a partial
self-covering $f:\M_1\arr\M$ is the quotient of the fundamental
group $\pi_1(\M, t)$ by the kernel of the iterated monodromy
action.
\end{definition}

\begin{figure}
 \includegraphics{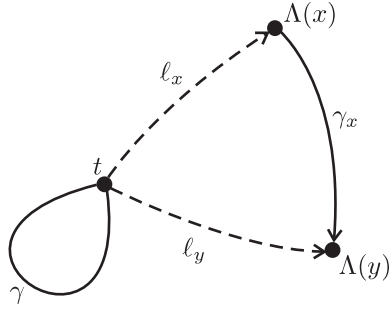}\\
 \caption{Self-similarity of iterated monodromy groups}\label{fig:recur}
\end{figure}

The iterated monodromy group is self-similar, if we identify the
tree of preimages $T$ with $\xs$ in a correct way
(see~\cite{nek:book,bgn}. The obtained self-similar action (called
the \emph{standard action}) is described in the following way.
Choose a bijection $\Lambda:\alb\arr f^{-1}(t)$ between the
alphabet $\alb$ of size $d$ and the set of preimages of the
basepoint. Choose also paths $\ell_x$ from the basepoint $t$ to
its preimage $\Lambda(x)$ for every $x$. Then the standard action
of $\img{f}$ on $\xs$ (which is conjugate to its action on $T$) is
given by the formula
\begin{equation}
\label{eq:imgstandard} \gamma(xw)=y(\ell_x\gamma_x\ell_y^{-1})(w),
\end{equation}
where $\gamma_x$ is the $f$-preimage of $\gamma$ starting in
$\Lambda(x)$, $y\in\alb$ is such that $\Lambda(y)$ is the end of
$\gamma_x$ and $w\in\xs$ is any word. We multiply here the paths
in the natural order (see Figure~\ref{fig:recur}).

\section{Self-similar algebras}

\subsection{Self-similarities of Hilbert spaces}

A ($d$-\emph{fold}) \emph{similarity} of an infinite-dimensional
Hilbert space $H$ is an isomorphism
\[\psi:H\arr H^d=\underbrace{H\oplus\cdots\oplus H}_{d}.\]

\begin{example}
Let $\alb$ be an alphabet with $d$ letters and let $\xs$ and $\xo$
be the rooted tree and its boundary, respectively. Let $\nu$ be
the uniform Bernoulli measure on $\xo$, i.e., the direct product
of uniform probability measures on $\alb$. Then the Hilbert space
$H=L^2(\xo, \nu)$ is decomposed into the direct sum
$\bigoplus_{x\in\alb}L^2(x\xo)$, where $L^2(x\xo)$ is the subspace
of functions with support a subset of the cylindrical set $x\xo$.
But the spaces $L^2(x\xo)$ are naturally isomorphic to $L^2(\xo,
\nu)$, where the isomorphism is the map $U_x:L^2(x\xo)\arr
L^2(\xo, \nu)$ given by
\[U_x(f)(w)=\frac 1{\sqrt{d}} f(xw).\] We view $U_x$ as
partial isometries of $L^2(\xo, \nu)$. Hence we get the natural
similarity $\sum_{x\in\alb} U_x$ of the space $L^2(\xo, \nu)$.
\end{example}

\begin{example}
Let $W\subset\xo$ be a \emph{self-similar} subset of the boundary
of the tree $\xs$, i.e., such a set that $W=\bigcup_{x\in\alb}xW$.
Then the space $\ell^2(W)$ is naturally self-similar, since it can
be decomposed into a direct sum
\[\ell^2(W)=\bigoplus_{x\in\alb}\ell^2(xW),\]
where the spaces $\ell^2(xW)$ are naturally isomorphic to
$\ell^2(W)$, with the isomorphism $U_x:\ell^2(xW)\arr\ell^2(W)$
given by
\[U_x(f)(w)=f(xw).\]
\end{example}

The above two examples will be the main types of self-similarities
on Hilbert spaces used in this paper.

\subsection{Representations of the Cuntz algebra $O_d$}

Recall that the \emph{Cuntz algebra} $O_d$, $d\ge 2$, is the
$C^*$-algebra given by the presentation
\[\langle a_1, a_2, \ldots, a_d\;:\;a_1a_1^*+a_2a_2^*+\cdots+a_da_d^*=1,\;a_k^*a_k=1,
k=1,\ldots, d\rangle.\] Note that as a corollary of the defining
relations we get $a_k^*a_l=0$ for $k\ne l$. The defining relations
can be also written as the following two matrix equalities
\[(a_1, a_2, \ldots, a_d)(a_1, a_2, \ldots, a_d)^*=1,\quad
(a_1, a_2, \ldots, a_d)^*(a_1, a_2, \ldots, a_d)=I,\] where $I$ is
the $d\times d$ unit matrix.

The Cuntz algebra $O_d$ is simple (see~\cite{cuntz}) and hence any
$d$ isometries $a_1, \ldots, a_d$ such that $a_1a_1^*+\cdots
+a_da_d^*=1$ generate a $C^*$-algebra isomorphic to the Cuntz
algebra $O_d$ and determine its representation.

Representations of the Cuntz algebra can be identified with
self-similarities of a Hilbert space as the following proposition
shows.

\begin{proposition}
\label{pr:cuntzrepr} The relation putting into correspondence to a
$*$-representation $\rho:O_d\arr B(H)$ the map
\[\tau_\rho=(\rho(a_1^*), \rho(a_2^*),
\ldots, \rho(a_d^*)):H\arr H^d\] is a bijection between the set of
representations of $O_d$ on $H$ and the set of $d$-fold
self-similarities on $H$.

The inverse of this bijection puts into correspondence to a
$d$-similarity $\psi:H\arr H^d$ the representation of $O_d$ given
by $\rho(a_k)=T_k$, for
\[T_k(\xi)=\psi^{-1}(0, \ldots, 0, \xi, 0, \ldots, 0),\]
where $\xi$ in the right-hand side is at the $k$th coordinate of
$H^d$.
\end{proposition}

\begin{proof}
If $\rho$ is a representation of $O_d$ on $H$, then $\rho(a_k)$
are isometries of $H$ with the subspaces $H_k=\rho(a_k)(H)$ and
$H=H_1\oplus \cdots \oplus H_d$. The $k$th components of a vector
$\xi\in H$ with respect to this decomposition is its image under
the projection $\rho(a_ka_k^*)$. Hence we get an isomorphism
\[H\arr H^d:\xi\mapsto (\rho(a_1)^*(\xi), \rho(a_2)^*(\xi), \ldots,
\rho(a_d)^*(\xi))\] equal to the composition of the identical map
$H\arr H_1\oplus\cdots\oplus H_d$ with the isomorphism
$\rho(a_1)^*\oplus\cdots\oplus\rho(a_d)^*:H_1\oplus\cdots\oplus
H_d\arr H^d$.

Conversely, suppose that
\[\psi:H\arr H^d\]
is a $d$-similarity. The map
\[\xi\mapsto(0, \ldots, 0, \xi, 0, \ldots 0)\]
is an isometry of $H$ with the $k$th direct summand of $H^d$.
Composing it with the isomorphism $\psi^{-1}$ we get an isometry
$T_k$ of $H$ with the direct summand $H_k$ of a decomposition
$H=H_1\oplus\cdots\oplus H_d$. But tuples of such isometries are
precisely the representations of the Cuntz algebra $O_d$. Note
that if $\psi(\xi)=(\xi_1, \ldots, \xi_d)$, then $\xi=\sum_k
T_k(\xi_k)$, hence $\xi_k=T_k^*(\xi)$. Consequently, the two
bijections are mutually inverse.
\end{proof}

\begin{example}\label{ex:murepr}
The representation of $O_d$ associated with the natural
$d$-similarity of $L^2(\xo, \nu)$ is generated by the isometries
$\pi(a_x)=T_x$ on $L^2(\xo, \nu)$, given by
\[T_x(f)(w)=\left\{\begin{array}{ll}\sqrt{d}f(w') &\text{if $w=xw'$}\\
0 & \text{otherwise.}\end{array}\right.\]
\end{example}

\begin{example}
Let $W$ be a self-similar subset of $\xo$. The representation
associated to the natural $d$-similarity on $\ell^2(W)$ is
generated by the isometries
\[T_x(f)(w)=\left\{\begin{array}{ll}f(w') &\text{if $w=xw'$}\\
0 & \text{otherwise.}\end{array}\right.\] Such representations of
$O_d$ are called \emph{permutational}. Permutational
representations of the Cuntz algebra related to self-affine
(digit) tilings of the Euclidean space are studied
in~\cite{cuntz_rep}.
\end{example}

\subsection{Self-similar groups and their representations}
\label{ss:wrrecurs} If $(G, \alb)$ is a self-similar group, then
the associated \emph{wreath recursion} is the embedding
$\phi:G\arr\symm\wr_\alb G=\symm\ltimes G^\alb$ given by
\[\phi(g)=\sigma(g|_x)_{x\in\alb},\]
where $\sigma\in\symm$ is the action $x\mapsto g(x)$ of $g$ on the
first level $\alb$ of the tree $\xs$ and the components $g|_x$ of
$G^\alb$ are the restrictions of $g$ onto the subtrees $x\xs$,
i.e., are given by the condition
\[g(xw)=g(x)g|_x(w)\]
for all $w\in\xs$.

Suppose that we have a self-similar group $(G, \alb)$ and let
$d=|\alb|$. Let $H$ be a Hilbert space together with a
$d$-similarity
\[\psi:H\arr H^\alb.\]
We will denote the summand of $H^\alb$ corresponding to a letter
$x\in\alb$ by $H_x$. Let $\rho:O_d\arr B(H)$ be the associated
representation of the Cuntz algebra. It is generated by isometries
$T_x=\rho(a_x)$, $x\in\alb$, such that $H_x=T_x(H)$.

\begin{definition}
A unitary representation $\rho$ of $G$ on $H$ is said to be
\emph{self-similar} (with respect to the $d$-similarity $\psi$) if
\[\rho(g)T_x=T_y\rho(h)\]
whenever $g(xw)=yh(w)$ for all $w\in\xs$, i.e., whenever $g\cdot
x=y\cdot h$ in the associated bimodule $\alb\cdot G$.
\end{definition}

\begin{example}
A self-similar group $G$ acts on the tree $\xs$ by automorphisms
and the induced action on the boundary $\xo$ is preserving the
Bernoulli measure $\nu$. We get hence a unitary representation
$\pi$ of $G$ on $L^2(\xo, \nu)$. It is easy to see that this
representation is self-similar with respect to the natural
$d$-similarity on $L^2(\xo, \nu)$.
\end{example}

\begin{example}
Let $W$ be a self-similar $G$-invariant subset of $\xo$. Then the
permutational representation of $G$ on $W$ is self-similar.
\end{example}

For the general notion of a Cuntz-Pimsner algebra
see~\cite{pimsner}.

\begin{definition}
Let $G$ be a self-similar group acing on $\xs$. The
\emph{associated (universal) Cuntz-Pimsner algebra} $\cpg$ is the
universal $C^*$-algebra generated by $G$ and $a_x$, $x\in\alb$,
satisfying the following relations
\begin{itemize}
\item all relations of $G$;
\item Cuntz relations for $a_x$: $a_x^*a_x=1$ for all $x\in\alb$,
$\sum_{x\in\alb}a_xa_x^*=1$;
\item $ga_x=a_yh$ for $g, h\in G$ and $x, y\in\alb$,
if $g(xw)=yh(w)$ for all $w\in\xs$, i.e., if $g(x)=y$ and
$h=g|_x$.
\end{itemize}

The algebra generated by $G$ in $\cpg$ is denoted $\maxalg$.
\end{definition}

Note that as a corollary of the defining relations we get the
relations
\begin{equation}\label{eq:splitting}g=g\sum_{x\in\alb}a_xa_x^*=
\sum_{x\in\alb}a_{g(x)}g|_xa_x^*,\end{equation} for
every $g\in G$, where $g|_x$ is, as usual, the section of $g$ at
$x$, i.e., such an element of $G$ that $g(xw)=g(x)g|_x(w)$ for all
$w\in\xs$.

The next proposition follows directly from the definitions.

\begin{proposition}
\label{pr:cpimsner} Let $\pi:O_d=\langle a_x\rangle_{x\in\alb}\arr
B(H)$ be a representation of the Cuntz algebra associated with a
$d$-similarity of a Hilbert space $H$.

A unitary representation $\rho$ of $G$ on $H$ is self-similar if
and only if $\rho$ and $\pi$ generate a representation of the
Cuntz-Pimsner algebra $\cpg$.

Consequently, self-similar representations of $G$ are precisely
restrictions onto $G$ of representations of the Cuntz-Pimsner
algebra $\cpg$.
\end{proposition}

\subsection{Matrix recursions}
A \emph{matrix recursion} on an algebra $A$ is a homomorphism
\[\phi:A\arr M_{d\times d}(A)\]
of $A$ into the algebra of matrices over $A$.

\begin{example}
Let $\psi:H\arr H^d=H^\alb$ be a $d$-similarity on a Hilbert space
$H$ and let $\rho:G\arr B(H)$ be a self-similar unitary
representation of a group $G$. Then every operator $\rho(g)$ for
$g\in G$ can be written, with respect to the decomposition
$\psi(H)=H^\alb$, as a $d\times d$ matrix
\[\rho(g)=\left(A_{yx}\right)_{x, y\in\alb},\]
where \[A_{yx}=\left\{\begin{array}{ll}\rho(g|_x) & \text{if $g(x)=y$,}\\
0 & \text{otherwise.}\end{array}\right.\]
\end{example}

For every self-similar group $G$ we have the associated matrix
recursion on the group algebra $\C[G]$, which is the linear
extension of the recursion:
\begin{equation}\label{eq:matrrecursion}\phi(g)=\left(A_{yx}\right)_{x,y\in\alb},\quad A_{yx}=\left\{\begin{array}{ll}g|_x & \text{if $g(x)=y$,}\\
0 & \text{otherwise,}\end{array}\right.\end{equation} which can be
interpreted as the \emph{wreath recursion} $\phi:G\arr
\symm\wr_\alb G$, as it was defined in Section~\ref{ss:wrrecurs}.

In terms of the associated representation $\rho$ of the Cuntz
algebra, we have
\[g|_x=T_{g(x)}^*gT_x,\]
where $T_x=\rho(a_x)$, which follows from~\eqref{eq:splitting} (or
directly from the defining relations of $\cpg$).

Note that the homomorphism $\phi$ usually is not injective (even
if the wreath recursion is and thus the map $\phi:G\arr
M_d\left(\C[G]\right)$ is injective as well).

\begin{example}
Consider the group $\mathfrak{G}=\langle a, b, c, d\rangle$ from
Example~\ref{sss:grigorchuk}. Then
\[\phi(a)=\left(\begin{array}{cc} 0 & 1 \\ 1 & 0\end{array}\right),
\phi(b)=\left(\begin{array}{cc} a & 0 \\ 0 & c\end{array}\right),
\phi(c)=\left(\begin{array}{cc} a & 0 \\ 0
& d\end{array}\right), \phi(d)=\left(\begin{array}{cc} 1 & 0 \\
0 & b\end{array}\right).\]

Denote $\alpha=(b+c+d-1)/2$. Then $\phi(\alpha)=\left(\begin{array}{cc} a & 0\\
0 & \alpha\end{array}\right)$,
$\phi(\alpha^2-1)=\left(\begin{array}{cc} 0 & 0\\ 0 &
\alpha^2-1\end{array}\right)$,
$\phi(a(\alpha^2-1)a)=\left(\begin{array}{cc} \alpha^2-1 & 0\\
0 & 0\end{array}\right)$, and
\[\phi\left((\alpha^2-1)a(\alpha^2-1)a\right)=0.\]

But
\[(\alpha^2-1)a(\alpha^2-1)a=\frac{(b+c+d-4)a(b+c+d-4)a}{16}\ne 0\]
in $\C[G]$.
\end{example}

One can find however an ideal $I\subset \C[G]$ such that $\phi$
induces an injective homomorphism $\C[G]/I\arr M_{d\times
d}(\C[G]/I)$, see~\cite{sid:ring,nek:iterat}. The ideal $I$ is the
ascending union of the kernels of the matrix recursions $\C[G]\arr
M_{d^n\times d^n}(\C[G])$ describing the action of the group $G$
on the $n$th level of the tree $\xs$.

\subsection{Self-similar completions of the group algebra}
Let $\rho:G\arr B(H)$ be a self-similar representation (with
respect to a $d$-similarity $\psi:H\arr H^d$) and let
$\mathcal{A}_\rho$ be the completion of $\C[G]$ with respect to
the norm given by $\rho$. Then the matrix recursion $\phi$ extends
to a homomorphism
\[\mathcal{A}_\rho\arr M_{d\times d}(\A_\rho),\]
also denoted by $\phi$ (or $\phi_\rho$), which is injective, since
it implements the equivalence of the representation $\rho$ with
the representation $\psi\circ\rho\circ\psi^{-1}$.

The following description of the completions $\mathcal{A}_\rho$
follows directly from Proposition~\ref{pr:cpimsner}.

\begin{definition}
A completion of $\C[G]$ is called \emph{self-similar} if it is the
completion with respect to a self-similar representation.
\end{definition}

The following proposition shows that there is a unique maximal
completion. We denote by $\maxalg$ the $C^*$-algebra generated by
$G$ in $\cpg$.

\begin{proposition}
\label{pr:algmaximal} A completion $\mathcal{A}$ of $\C[G]$ is
self-similar if and only if it is the closure of $\C[G]$ in a
homomorphic image of the Cuntz-Pimsner algebra $\cpg$. In
particular, every such completion is a homomorphic image of the
algebra $\maxalg$.
\end{proposition}

 There also exists the smallest self-similar completion.

\begin{definition}\label{def:Ggeneric}
Let $G$ be a countable group of automorphisms of $\xs$. A point
$w\in\xo$ of the boundary is called \emph{$G$-generic} if for
every $g\in G$ either $g(w)\ne w$ or $w$ is fixed by $g$ together
with all points of a neighborhood of $w$.
\end{definition}

It is not hard to prove that the set of $G$-generic points is
co-meager (i.e., is an intersection of a countable collection of
open dense sets), see~\cite{grineksu_en} and~\cite{nek:bim}. Let
$W\subset\xo$ be a non-empty countable $G$-invariant set of
$G$-generic points (one can also just take the $G$-orbit of a
$G$-generic point). Denote by $\rho_W$ the permutational
representation of $G$ on $\ell^2(W)$. The following theorem is
proved in~\cite{nek:bim}.

\begin{theorem}\label{th:minimal} Let $G$ be a group with a self-similar action on
$\xs$. Denote by $\|\cdot\|_{min}$ the norm on $\C[G]$ defined by
the representation $\rho_W$. This norm does not depend on the
choice of the set $W$ and $\|a\|_{min}\le\|a\|_{\rho}$ for any
self-similar representation $\rho$ of $G$ and any $a\in\C[G]$.

The completion of $\C[G]$ with respect to the norm
$\|\cdot\|_{min}$ is the algebra $\minalg$ generated by $G$ in a
unique simple unital quotient of $\cpg$.
\end{theorem}

Hence, if $\mathcal{A}$ is a completion of $\C[G]$ with respect to
a self-similar representation of $G$, then the identical map on
$G$ induces surjective homomorphisms
\[\maxalg\arr\mathcal{A}\arr\minalg.\]

Let $\mathfrak{O}_G$ be the groupoid generated by the germs of the
local homeomorphisms of $\xo$ of the form
\begin{equation}\label{eq:tvg}T_vg:w\mapsto vg(w),\quad v\in\xs, g\in G.\end{equation}
Recall that a germ of a local homeomorphism $h:U\arr V$ is a pair
$(h, x)$, where $x$ belongs to the domain of $h$, where two pairs
$(h_1, x_1)$ and $(h_2, x_2)$ are identified if $x_1=x_2$ and
restrictions of $h_i$ on a neighborhood of $x_1$ coincide. The
germs are composed in a natural way and inverse of a germ $(g, x)$
is defined to be equal to $(g^{-1}, g(x))$. The set of germs of a
pseudogroup of local homeomorphisms has a natural \emph{germ
topology} defined by the basis consisting of the sets of the form
$\mathcal{U}_{h, U}=\{(h, x)\;:\;x\in U\}$ where $h$ is an element
of the pseudogroup and $U$ is an open subset of the domain of $h$.

For more details and for the definition of the operator algebras
associated with topological groupoids,
see~\cite{renault:groupoids,bridhaefl,paterson:gr}.

It is not hard to prove that the universal Cuntz-Pimsner algebra
$\cpg$ coincides with the universal algebra of the groupoid
$\mathfrak{O}_G$ (see~\cite{nek:cpalg}). One can consider also the
\emph{reduced $C^*$-algebra} of $\mathfrak{O}_G$ (see the
definitions in~\cite{renault:groupoids,khoshkamskand}). Let
$\A_{red}$ be the subalgebra of the reduced $C^*$-algebra of
$\mathfrak{O}_G$ generated by $G$. Then $\A_{red}$ is also a
completion of $G$ with respect to a self-similar representation
(since it comes from a representation of $\cpg$, see
Proposition~\ref{pr:cpimsner}).

A groupoid of germs of an action of a group $G$ on a topological
space $\mathcal{X}$ is non-Hausdorff if and only if there exists
$x\in\mathcal{X}$ and an element $g\in G$ such that the germ of
$g$ at $x$ can not be separated from the germ of the identity at
$x$. The latter is equivalent to the condition that for every
neighborhood $U$ of $x$ there exists $y\in U$ such that the germ
of $g$ at $y$ is trivial (in particular $g(x)=x$) and there exists
$z\in U$ such that $g(z)\ne z$. Consequently, the groupoid of
germs is Hausdorff if and only if for every $g\in G$ the interior
of the set of fixed points of $g$ is closed. See an example of a
self-similar group with non-Hausdorff groupoid of germs of the
action on $\xo$ in Example~\ref{ex:nonHausd}.

The following theorem follows from classical results of theory of
amenable groupoids
(see~\cite{renault:groupoids,delaroche_renault}). The last
paragraph of the theorem follows directly from the definitions of
the corresponding algebras.

\begin{theorem}
\label{th:amenablegr} If the groupoid of germs of the action of
$G$ on $\xo$ is (measurewise) amenable (in particular, if the
orbits of the action of $G$ on $\xo$ have polynomial growth), then
the universal and reduced algebras of the groupoid
$\mathfrak{O}_G$ coincide. In particular, their sub-algebras
$\maxalg$ and $\A_{red}$ coincide.

If the groupoid of germs of the action of $G$ on $\xo$ is
Hausdorff, then the algebras $\A_{red}$ and $\minalg$ coincide.
\end{theorem}

In particular, if the action of $G$ on $\xo$ is free, then the
groupoid of germs is Hausdorff and moreover is \emph{principal}
(i.e., is an equivalence relation) and $\A_{red}=\minalg$.

Another frequently used self-similar completion is the completion
$\A_{mes}$ defined by the natural unitary representation $\pi$ of
$G$ on $L^2(\xo, \nu)$.

\begin{question} Find conditions when $\A_{mes}=\minalg$.\end{question}

\subsection{Gauge-invariant subalgebra of $\cpg$}
Let $\M_k$ be the closed linear span in $\cpg$ of the elements
$a_vga_u^*$ for $g\in G$ and $v, u\in\alb^k$. Here we use the
multi-index notation $a_{x_1x_2\ldots x_n}=a_{x_1}a_{x_2}\cdots
a_{x_n}$. In particular, $\M_0=\maxalg$ is the algebra generated
by $G$ in $\cpg$.

It is easy to see that for $v_1, v_2, u_1, u_2\in\alb^k$ and $g_1,
g_2\in G$
\[a_{v_1}g_1a_{u_1}^*a_{v_2}g_2a_{u_2}^*=\left\{\begin{array}{cr}a_{v_1}g_1g_2a_{u_2}^*
& \text{if $u_1=v_2$}\\
0 & \text{otherwise.}\end{array}\right.\] Observe that $a_va_u^*$
are multiplied in the same way as the matrix units, hence $\M_k$
is isomorphic to the algebra $M_{d^k\times d^k}(\maxalg)$ of
$d^k\times d^k$-matrices over $\maxalg$.

Recall that by~\eqref{eq:splitting} every element $g\in G$ can be
written in $\cpg$ as a sum
\[
g=\sum_{x\in\alb}a_{g(x)}g|_x a_x^*.
\]
If we apply this formula to the element $g$ in $a_vga_u^*\in\M_k$,
we get an element of $\M_{k+1}$. Hence $\M_{k+1}\supset\M_k$.

The algebra $\gi$ is defined as the closure in $\cpg$ of the
ascending union $\bigcup_{k\ge 0}\M_k$.

Every algebra $\M_k$ contains the sub-algebra equal to the linear
span of the elements $a_va_u^*$, which is isomorphic to
$M_{d^k\times d^k}(\C)$. Their union (i.e., the closed linear span
of all elements $a_va_u^*$ for $|v|=|u|$) is isomorphic to the
Glimm's uniformly hyperfinite algebra $M_{d^\infty}(\C)$
(see~\cite{glimm,davidson}), since the expansion
rule~\eqref{eq:splitting} for $g=1$ defines the diagonal
embeddings $M_{d^k\times d^k}\hookrightarrow M_{d^{k+1}\times
d^{k+1}}$. Its diagonal subalgebra, generated by the projections
$a_va_v^*$ is isomorphic to the algebra $C(\xo)$ of continuous
functions on $\xo$. The isomorphism identifies a projection
$a_va_v^*$ with the characteristic function of the cylindrical set
$v\xo$. It is easy to see that the action of $G\subset\gi$ on the
diagonal algebra by conjugation coincides with the action by
conjugation of $G$ on $C(\xo)$.

It is proved in~\cite{nek:bim} that if the group $(G, \alb)$ is
self-replicating, then the subalgebra $\gi$ of $\cpg$ is generated
by $G$ and the subalgebra $C(\xo)$, i.e., that $\gi$ is the
cross-product of $\maxalg$ and the algebra of continuous functions
on $\xo$ induced by the usual action of $G$ on $\xo$.

Similarly to the Cuntz algebra (see~\cite{cuntz,davidson}), we
have a natural strongly continuous (\emph{gauge}) action $\Gamma$
of the circle $\mathbb{T}=\{z\in\C: |z|=1\}$ on $\cpg$ by
\begin{eqnarray*}
\Gamma_z(g)  & = & g \\
\Gamma_z(a_x) & = & za_x.
\end{eqnarray*}
for $g\in G, x\in\alb$ and $z\in\mathbb{T}$.

 Then $\Gamma_z(a_vga_u^*)=z^{|v|-|u|}a_vga_u^*$ for $u, v\in\xs$, $g\in G$,
thus the integral $\int\Gamma_z(a_vga_u^*)dz$ is equal to zero for
$|v|\neq |u|$ and to $a_vga_u^*$ for $|v|=|u|$, where $dz$ is the
normalized Lebesgue measure on the circle. Consequently the map
\begin{equation}
\label{eq:exp} \mathbb{M}_0(a)=\int\Gamma_z(a)dz
\end{equation}
for $a\in\cpg$ is a conditional expectation from $\cpg$ onto the
subalgebra $\gi$.

\begin{example}
If $G$ is the cyclic group generated by the adding machine
$a=\sigma(1, a)$, then the algebra $\maxalg=\minalg$ is isomorphic
to the algebra $C(\mathbb{T})$ with the linear recursion
$C(\mathbb{T})\arr M_2(C(\mathbb{T}))$ coming from the double
self-covering of the circle (if we identify $C(\mathbb{T})$ with
$C^*(\Z)$ via Fourier series, then this linear recursion is given
by $z\mapsto\left(\begin{array}{cc} 0 & z\\ 1 &
0\end{array}\right)$, where $z=e^{2\pi i t}$, $t\in\R$, is the
variable of the Fourier series). Consequently, the algebra $\gi$
of the adding machine action is the Bunce-Deddence algebra. It is
also isomorphic to the cross-product algebra of the odometer
action on the Cantor space $\xo$ (see~\cite{davidson}).
\end{example}

\subsection{Overview of algebras associated with self-similar groups}

We studied above the universal Cuntz-Pimsner algebra $\cpg$ of a
self-similar group $G$. Representations of this algebras
correspond to representations of $G$ on a Hilbert space $H$ which
are self-similar with respect to some similarity of $H$ (see
Proposition~\ref{pr:cpimsner}).

The subalgebra of $\cpg$ generated by $G$ was denoted $\maxalg$.
Any self-similar representation of $G$ extends to a representation
of $\maxalg$, hence $\maxalg$ is the completion of the group
algebra with respect to the maximal self-similar norm.

Third algebra is the natural inductive limit of the algebras
$M_{d^n\times d^n}(\maxalg)$, which we denoted by
$\M_{d^\infty}(G)$. It is a subalgebra of $\cpg$ in a natural way.

Hence we get the following tower of algebras associated to a
self-similar group
\[\maxalg\subset\M_{d^\infty}(G)\subset\cpg.\]

There exits also the smallest self-similar norm on the group
algebra, which corresponds to a unique simple homomorphic image
${\cpg}_{min}$ of the Cuntz-Pimsner algebra $\cpg$ (see
Theorem~\ref{th:minimal}). The subalgebra of ${\cpg}_{min}$
generated by $G$ is denoted $\minalg$ and it is the completion of
the group algebra with respect to the smallest self-similar norm.
The image $\M_{d^\infty}(G)_{min}$ of $\M_{d^\infty}(G)$ is also
simple.

For any other self-similar representation $\rho$ of $G$ we get the
corresponding homomorphic image of $\cpg$, $\M_{d^\infty}(G)$ and
$\maxalg$. We get hence the following commutative diagram of
algebras
\[\begin{array}{ccccc}
\maxalg & \hookrightarrow & \M_{d^\infty}(G) & \hookrightarrow & \cpg\\
\mapdown{} & & \mapdown{} & & \mapdown{} \\
\A_\rho & \hookrightarrow & \M_{d^\infty}(G)_\rho & \hookrightarrow & {\cpg}_\rho\\
\mapdown{} & & \mapdown{} & & \mapdown{} \\
\minalg & \hookrightarrow & \M_{d^\infty}(G)_{min} &
\hookrightarrow & {\cpg}_{min}
\end{array}\]
where all vertical arrows are surjective and all horizontal are
embeddings.

An important example of a self-similar representation is the
natural unitary representation $\pi$ of $G$ on $L^2(\xo, \nu)$,
where $\nu$ is the uniform Bernoulli measure on the boundary $\xo$
of the rooted tree. Together with the natural self-similarity of
the space $L^2(\xo, \nu)$ this gives a representation of $\cpg$.
We denote the respective homomorphic images by
$\A_{mes}\hookrightarrow
\M_{d^\infty}(G)_{mes}\hookrightarrow{\cpg}_{mes}$.

The algebra $\A_{mes}$ is residually finite-dimensional, since the
representation $\pi$ of $G$ is a direct sum of finite-dimensional
representations (coming from the action of $G$ on the levels of
the tree). These finite-dimensional representations give an
additional tool to the study of the algebra $\A_{mes}$.

Therefore the next questions about the epimorphisms
$\maxalg\arr\A_{mes}\arr\minalg$ are natural.

\begin{question}
When are the epimorphisms $\A_{mes}\arr\minalg$ and
$\maxalg\arr\A_{mes}$ isomorphisms?
\end{question}

\begin{question}
Under which conditions the algebra $\A_{mes}$ isomorphic to the
reduced $C^*$-algebra of the group $G$? Is it so when $G$ is the
free group generated by the automaton shown on
Figure~\ref{fig:alfree} in Example~\ref{sss:free}?
\end{question}

\section{Contracting groups and limit solenoids}
\label{s:contracting}
\subsection{Matrix recursion for contracting groups}

\begin{definition}
A self-similar group $(G, \alb)$ is said to be \emph{contracting}
if there exists a finite set $\nuke\subset G$ such that for every
$g\in G$ there exists $n_0$ such that $g|_v\in\nuke$ for all words
$v\in\xs$ of length greater than $n_0$. The smallest set $\nuke$
with this property is called the \emph{nucleus} of the
self-similar group.
\end{definition}

\begin{example}
The adding machine action of $\Z$ is contracting with the nucleus
$\nuke=\{1, a, a^{-1}\}$, since $a^n|_0=a^{\lfloor n/2\rfloor}$
and $a^n|_1=a^{\lceil n/2 \rceil}$.
\end{example}

\begin{example}
It is also not hard to prove that the torsion group from
Example~\ref{sss:grigorchuk} is contracting with the nucleus $\{1,
a, b, c, d\}$.
\end{example}

\begin{example}
The free group considered in Example~\ref{sss:free} is not
contracting. In this example any section $g|_v$ has the same
length as $g$.
\end{example}

The nucleus $\nuke$ of a contracting group can be interpreted as
an \emph{automaton}, i.e., for every $g\in\nuke$ and every letter
$x\in\alb$ the restriction $g|_x$ belongs to $\nuke$.
Consequently, if we denote by $N$ the linear span of $\nuke$ in
$\C[G]$ and $\phi:\C[G]\arr M_{d\times d}(\C[G])$ is the
associated matrix recursion, then $\phi(N)$ is a subspace of the
space $M_{d\times d}(N)$ of matrices with entries in $N$.
Moreover, the following is true.

\begin{theorem}[\cite{nek:cpalg}]
If the group $(G, \alb)$ is contracting and $\nuke$ is its
nucleus, then the algebra $\cpg$ is generated by
$\{a_x\}_{x\in\alb}\cup\nuke$ and is defined by the following
finite set of relations
\begin{enumerate}
\item Cuntz relations \[a_x^*a_x=1\]
\item decompositions
\[g=\sum_{x\in\alb}a_{g(x)}g|_xa_x^*\]
for every $g\in\nuke$ (this includes the remaining Cuntz algebra
relation \[\sum_{x\in\alb}a_xa_x^*=1\] in the case $g=1$),
\item all relations $g_1g_2g_3=1$ of length at most three which are true for the elements
of the nucleus $\nuke$ in the group $G$ and relations
$gg^*=g^*g=1$ for $g\in\nuke$.
\end{enumerate}
\end{theorem}

Moreover, the groupoid of germs of the action of $G$ on $\xo$ is
amenable in the case of a contracting group $G$. This follows from
the fact that the orbits of the of the action on $\xo$ have
polynomial growth (see~\cite{delaroche_renault}
and~\cite{nek:book}). This implies the following corollary of
Theorem~\ref{th:amenablegr}.

\begin{proposition}
If the group $G$ is contracting and level-transitive, then the
algebras $\maxalg$ and $\A_{red}$ coincide. If, additionally, for
every element $g$ of the nucleus of $G$ the interior of the set of
fixed points of $g$ is closed, then all self-similar completions
of $G$ are isomorphic to $\maxalg=\minalg=\A_{red}$.
\end{proposition}

\begin{examples}
The adding machine action, the Basilica group~\ref{sss:basilica}
satisfy all the conditions of the corollary, hence they have
unique self-similar completions.
\end{examples}

\begin{example}
\label{ex:nonHausd} Consider the group generated by the
transformations $a, b, c$ given by
\begin{alignat*}{2}
a(0w) &= 1w, &\quad a(1w) &= 0w\\
b(0w) &= 0a(w), &\quad b(1w) &=1c(w)\\
c(0w) &= 0w, &\quad c(1w) &= 1b(c).
\end{alignat*}
This is one of the Grigorchuk groups $G_w$ studied
in~\cite{grigorchuk:growth_en} (for $w=111\ldots$). Its growth was
studied by A.~Erschler in~\cite{ershler:growth}.

This group is contracting with the nucleus $\{1, a, b, c,
bc=cb\}$, but its groupoid of germs is not Hausdorff. Namely, the
set of fixed points of $c$ is equal to
\[\{111\ldots\}\cup\bigsqcup_{k=0, 1, \ldots}\underbrace{11\ldots
1}_{2k}0\xo,\] its interior is $\bigsqcup_{k=0, 1,
\ldots}\underbrace{11\ldots 1}_{2k}0\xo$, which is not closed.

Let us show that in this case the element $b+c-bc-1$ belongs to
the kernel of the epimorphism $\maxalg=\A_{red}\arr\minalg$.

We have \begin{multline*}\phi(b+c-bc-1)=\left(\begin{array}{cc} a & 0\\
0 & c\end{array}\right)+\left(\begin{array}{cc} 1 & 0\\ 0 &
b\end{array}\right)-\left(\begin{array}{cc} a & 0\\ 0 &
bc\end{array}\right)-\left(\begin{array}{cc} 1 & 0\\ 0 &
1\end{array}\right)=\\ \left(\begin{array}{cc} 0 & 0\\ 0 &
b+c-bc-1\end{array}\right),\end{multline*} which implies that
$\pi(b+c-bc-1)=0$ for the permutation representation $\pi$ on any
orbit of $G$ on $\xo$. Hence $b+c-bc-1$ is equal to zero in
$\minalg$.

On the other hand, it follows from contraction that the isotropy
group of the point $111\ldots\in\xo$ in the groupoid of germs of
the action of $G$ on $\xo$ contains 4 elements (the germs of $1,
b, c$ and $bc$). Consequently, for any germ $\gamma$ with range
$111\ldots$ the germs $\gamma$, $b\cdot\gamma$, $c\cdot\gamma$ and
$bc\cdot\gamma$ are pairwise different (here the \emph{range} of a
germ $(h, x)$ is the point $h(x)$). This implies that the element
$b+c-bc-1$ is non-zero in the regular representation of the
groupoid, hence it is non-zero in $\A_{red}$.

Consequently, the homomorphism $\A_{red}\arr\minalg$ is not an
isomorphism in general. In particular, $\A_{red}$ is not simple,
even though the action of the group $G$ is minimal. (If the action
is minimal and the groupoid of germs is Hausdorff, then the
reduced algebra of the groupoid of germs is simple,
see~\cite{renault:groupoids} Proposition~4.6.)
\end{example}

\begin{question} Describe the kernel of the epimorphism
$\A_{red}\arr\minalg$. Is it true that for a contracting group $G$
it is generated (in $\cpg$) by linear combinations of the elements
of the nucleus?
\end{question}

\subsection{Limit solenoid}
Let us fix some contracting self-similar group $(G, \alb)$ with
the nucleus $\nuke$. Consider the space $\xz$ of bi-infinite
sequences of the form
\[\ldots x_{-2}x_{-1}\;.\;x_0x_1\ldots\]
of letters $x_i\in\alb$. Here the dot marks the place between the
coordinates number 0 and number $-1$. We consider $\xz$ to be a
topological space with the direct product topology of discrete
sets $\alb$.

\begin{definition}
Two sequences $\ldots x_{-2}x_{-1}\;.\;x_0x_1\ldots$ and $\ldots
y_{-2}y_{-1}\;.\;y_0y_1\ldots$ are said to be \emph{asymptotically
equivalent} (with respect to the action of $G$) if there exists a
finite set $N\subset G$ and a sequence $g_i\in N$ such that
\[g_i(x_ix_{i+1}x_{i+2}\ldots)=y_iy_{i+1}y_{i+2}\ldots\]
for all $i\in\mathbb{Z}$.
\end{definition}

It is proved in~\cite{nek:book} that we can take $N$ equal to the
nucleus of $G$ and that the asymptotic equivalence relation can be
described in the following way.

\begin{proposition}
\label{pr:aseq} The sequences $\ldots
x_{-2}x_{-1}\;.\;x_0x_1\ldots$ and $\ldots
y_{-2}y_{-1}\;.\;y_0y_1\ldots$ are asymptotically equivalent if
and only if there exists a sequence $g_i\in\nuke$ of elements of
the nucleus such that $g_i\cdot x_i=y_i\cdot g_{i-1}$, i.e., such
that $y_i=g_i(x_i)$ and $g_{i-1}=g_i|_{x_i}$.
\end{proposition}

\begin{definition}
The \emph{limit solenoid} $\solen$ of the self-similar group $(G,
\alb)$ is the quotient of the topological space $\xz$ by the
asymptotic equivalence relation.
\end{definition}

The next proposition follows from the definition and
Proposition~\ref{pr:aseq} (see details in~\cite{nek:book}).

\begin{proposition}
The limit solenoid $\solen$ is a compact metrizable
finite-dimensional space. If the action of $G$ on $\xs$ is
level-transitive, then $\solen$ is connected. The shift
\[\ldots x_{-2}x_{-1}\;.\;x_0x_1\ldots\mapsto
\ldots x_{-3}x_{-2}\;.\;x_{-1}x_0\ldots\] induces a homeomorphism
$\si:\solen\arr\solen$.
\end{proposition}

\begin{definition}
Let $\xmo$ be the space of sequences $\ldots x_2x_1$, over the
alphabet $\alb$ with the direct product topology. Two sequences
$\ldots x_2x_1, \ldots y_2y_1\in\xmo$ are \emph{asymptotically
equivalent} if there exists a finite set $N\subset G$ and a
sequences $g_k\in N$ such that $g_k(x_k\ldots x_1)=y_k\ldots y_1$
for all $k\ge 1$. The quotient of $\xmo$ by the asymptotic
equivalence relation is called the \emph{limit space} of the group
$(G, \alb)$ and is denoted $\lims$.
\end{definition}

Here also we can take $N$ to be equal to the nucleus
(see~\cite{nek:book}).

We have a natural continuous projection $\solen\arr\lims$ induced
by the map \[\ldots x_{-2}x_{-1}\;.\;x_0x_1\ldots\mapsto \ldots
x_{-2}x_{-1}.\] This projection semiconjugates the map
$\si:\solen\arr\solen$ with the map $\sis:\lims\arr\lims$ induced
by the one-sided shift $\ldots x_2x_1\mapsto \ldots x_3x_2$. We
call $(\lims, \sis)$ the \emph{limit dynamical system}.

The following theorem is proved in~\cite{nek:book}, where a more
general formulation can be found.

\begin{theorem}
If $f:\M_1\arr\M$ is an expanding partial covering (where $\M_1$
is an open subset of $\M$ with the induced Riemann metric), then
the iterated monodromy group $\img{f}$ is contracting and the
limit dynamical system $\sis:\lims[\img{f}]\arr\lims[\img{f}]$ is
topologically conjugate to the action of $f$ on its Julia set.
\end{theorem}

The limit solenoid $\solen$ can be reconstructed from the limit
dynamical system as the inverse limit of the sequence
\[\lims\stackrel{\sis}{\longleftarrow}\lims\stackrel{\sis}{\longleftarrow}\ldots.\]
The map $\sis$ induces a homeomorphism of the inverse limit, which
is conjugate with $\si:\solen\arr\solen$.

\begin{example}[Lyubich-Minsky laminations]
Let $f\in\C(z)$ be a post-critically finite rational function
(i.e., the orbits of its critical points are finite). Consider the
inverse limit $\widehat{\mathcal{S}_f}$ of the backward iteration
\[\widehat\C\stackrel{f}{\longleftarrow}\widehat\C\stackrel{f}{\longleftarrow}\ldots,\]
where $\widehat\C$ is the Riemann sphere. The shift along the
inverse sequence (i.e., along the action of $f$) defines a
homeomorphism $\widehat
f:\widehat{\mathcal{S}_f}\arr\widehat{\mathcal{S}_f}$, called the
\emph{natural extension} of $f$. We have the natural projection
$P$ of $\widehat{\mathcal{S}_f}$ onto $\widehat\C$ (onto the first
term of the inverse sequence).

Let $\mathcal{S}_f$ be the preimage of the Julia set in
$\widehat{\mathcal{S}_f}$ under the projection $P$. Then the space
$\mathcal{S}_f$ is homeomorphic to the limit solenoid of $\img{f}$
and the action of $\widehat f$ on $\mathcal{S}_f$ is topologically
conjugate to the action of the shift on the limit solenoid.

The space $\widehat{\mathcal{S}_f}$ and the homeomorphism
$\widehat f$ were studied in~\cite{lyubichminsk}.
\end{example}

\subsection{The solenoid as a hyperbolic dynamical system}
Let us have a look at the dynamical system $\si:\solen\arr\solen$
more carefully. To avoid technicalities, we will assume that the
self-similar group $(G, \alb)$ is finitely generated, contracting,
self-replicating and that it is \emph{regular}, which means that
for every $g\in G$ and $w\in\xo$ either $g(w)\ne w$, or $w$ is
fixed by $g$ together will all points of a neighborhood of $w$
(i.e., for some beginning $v$ of $w$ the restriction $g|_v$ is
trivial). In other words, a group $G$ is said to be regular if
every point of $\xo$ is $G$-regular in the sense of
Definition~\ref{def:Ggeneric}. It is sufficient to check the
regularity condition for the elements $g$ belonging to the nucleus
of $G$.

If $f:\mathcal{M}\arr\mathcal{M}$ is an expanding self-covering of
a complete compact geodesic space, then the iterated monodromy
group $G=\img{f}$ is contracting and the dynamical system
$\sis:\lims[\img{f}]\arr\lims[\img{f}]$ is topologically conjugate
to $(\mathcal{M}, f)$. Moreover, $\img{f}$ satisfies the
regularity condition in this case (this easily follows from the
description of the orbispace structure on $\lims[\img{f}]$ given
in~\cite{nek:book}).

We say that two points $\xi, \zeta\in\solen$ are \emph{stably
(unstably) equivalent} if for every neighborhood of the diagonal
$U\subset\solen\times\solen$ there exists $n_U\ge 0$ such that
\[(\si^n(\xi), \si^n(\zeta))\in U\] for all $n\ge n_U$ ($n\le -n_U$,
respectively).

\begin{proposition}
\label{pr:reghyperb} Let the points $\xi, \zeta\in\solen$ be
represented by sequences $(x_n)_{n\in\Z}$ and $(y_n)_{n\in\Z}$,
respectively.

The points $\xi, \zeta$ are stably equivalent if and only if there
exists $n\in\Z$ such that the sequences $\ldots x_{n-1}x_n$ and
$\ldots y_{n-1}y_n\in\xmo$ are asymptotically equivalent, i.e.,
represent the same point of $\lims$.

The points $\xi, \zeta$ are unstably equivalent if and only if
there exists $n\in\Z$ such that
\[g(x_nx_{n+1}\ldots)=y_ny_{n+1}\ldots\] for some element $g$ of the nucleus.
\end{proposition}

\begin{proof}
It is easy to see that stable and unstable equivalence follows
from the conditions in the proposition.

Let us show that the converse implications hold. Let $k$ be such
that for every two elements $g, h\in\nuke$ of the nucleus and
every word $v\in\alb^n$ either $g(v)\ne h(v)$ or $g(v)=h(v)$ and
$g|_v=h|_v$. Such $k$ exists by the regularity condition.

Denote by $\overline U_k$ the set of pairs $\left((x_n)_{n\in\Z},
(y_n)_{n\in\Z}\right)$ such that $g(x_{-k}\ldots x_k)=y_{-k}\ldots
y_k$ for some element $g\in\nuke$ of the nucleus. Let $U_k$ be the
image of $\overline U_k$ in $\solen\times\solen$.

Every set $U_k$ is a neighborhood of the diagonal and for every
neighborhood of the diagonal $U\subset\solen\times\solen$ there
exists $k$ such that $U_k\subset U$.

Let $\xi$ and $\zeta$ be the points of $\solen$ represented by the
sequences $(x_i)_{i\in\Z}$ and $(y_i)_{i\in\Z}$.

Suppose that $\xi$ and $\zeta$ are stably equivalent. Then for any
$k$ there exists $n_k$ such that the pair
$\left(\si^n\left((x_i)_{i\in\Z}\right),
\si^n\left((y_i)_{i\in\Z}\right)\right)$ belongs to $\overline
U_k$ for all $n\ge n_k$. Then for every $n>n_k+k$ there exists
$g_n\in\nuke$ such that $g_n(x_{-n}x_{-n+1}\ldots
x_{-n+2k})=y_{-n}y_{-n+1}\ldots y_{-n+2k}$. Consequently, for
$h_n=g_n|_{x_{-n}\ldots x_{-n+k-1}}$ we have \[h_n(x_{-n+k}\ldots
x_{-n+2k})=y_{-n+k}\ldots y_{-n+2k}.\] Note that by the choice of
$k$ the element $h_n$ is determined uniquely by the pair of words
$(x_{-n}x_{-n+1}\ldots x_{-n+2k}, y_{-n}y_{-n+1}\ldots y_{-n+2k})$
(i.e., does not depend on the choice of $g_n$). The uniqueness of
$g_n$ implies that $h_n\cdot x_{-n+k}=y_{-n+k}\cdot h_{n-1}$,
which finishes the proof.

The proof of the statement about the unstable equivalence is
analogous.
\end{proof}

Hence we get the following.

\begin{corollary}
In conditions of Proposition~\ref{pr:reghyperb} the points $\xi$
and $\zeta$ are unstably equivalent if and only if $x_0x_1\ldots$
and $y_0y_1\ldots$ belong to one $G$-orbit.

In other words, the unstable equivalence relation
$\mathcal{U}\subset\solen\times\solen$ is a union of the sets
$\mathcal{U}_g$ equal to the images in $\solen\times\solen$ of the
sets
\[\{(\ldots x_{-2}x_{-1}\;.\;x_0x_1\ldots; \ldots
y_{-2}y_{-1}\;.\;g(x_0x_1\ldots))\;:\;x_i, y_i\in\alb\}.\]
\end{corollary}

Let us introduce a topology on $\mathcal{U}$ equal to the (direct
limit) topology of the union $\mathcal{U}=\bigcup_{g\in
G}\mathcal{U}_g$, where $\mathcal{U}_g$ are taken with the induced
topology of the subset $\mathcal{U}_g\subset\solen\times\solen$
and the direct limit is taken with respect to the identical maps
between the sets $\bigcup_{g\in A}\mathcal{U_g}$, where $A$ runs
through the finite subsets of $G$.

This topology converts the unstable equivalence relation into a
Hausdorff topological groupoid, which we call the \emph{unstable
groupoid}. A natural Haar system on it comes from the push forward
of the uniform Bernoulli measure $\nu$ on $\xmo$.

The shift $\si:\solen\arr\solen$ defines automorphisms of the
stable and unstable groupoids.

The next theorem can be seen as a formulation of the fact that $G$
is the holonomy group of the unstable foliation of $\solen$ in
terms of $C^*$-algebras.

\begin{theorem}
\label{th:moritaequiv} The universal convolution algebra of the
unstable groupoid is strongly Morita equivalent to $\gi$.
\end{theorem}

\begin{proof}
Let us construct an equivalence bimodule (in the sense
of~\cite{muhlyrenault:equiv}) between the groupoid $\mathcal{U}$
of the unstable equivalence relation and the groupoid
$\mathcal{G}$ of germs of the semigroup of transformations
$T_vgT_u^*$ for $u, v\in\xs, |u|=|v|$ and $g\in G$. The universal
convolution algebra of $\mathcal{G}$ is $\gi$
(see~\cite{nek:cpalg}). Note also that regularity of $(G, \alb)$
implies that the groupoid $\mathcal{G}$ is principal (i.e., is an
equivalence relation), since the isotropy groups of $\mathcal{G}$
are trivial.

Let us consider the bimodule $Z$ consisting of the set of pairs
$(\zeta, w)$, where $\zeta\in\solen$ and $w\in\xo$ are such that
$\ldots x_{-2}x_{-1}\;.\;w$ represents a point unstably equivalent
to $\zeta$ (for some and hence for any $\ldots x_{-2}x_{-1}$). Let
$Z_g\subset Z$ for $g\in G$ be the set of pairs $(\zeta, w)$ such
that $\zeta$ is represented by a sequence $\ldots
x_{-2}x_{-1}\;.\;g(w)$. Endow $Z_g$ with the induced topology from
$\solen\times\xo$. The set $Z$ is a union of $Z_g$ for all $g\in
G$. Consider the direct limit topology on $Z$ coming from this
union.

We have a natural right action of $\mathcal{U}$ on $Z$:
\[(\zeta_2, w)\cdot (\zeta_1, \zeta_2)=(\zeta_1, w)\]
and a left action of $\mathcal{G}$ on $Z$:
\[(w_2, w_1)\cdot (\zeta, w_1)=(\zeta, w_2),\]
which satisfy the conditions of the equivalence bimodule (see
Definition~2.1 of~\cite{muhlyrenault:equiv}).

One can find another proof of this theorem in~\cite{nek:cpalg}.
\end{proof}

The dynamical system $\si:\solen\arr\solen$ is an example of a
\emph{Smale space} (for a definition see~\cite{put_kam_sp}). The
$C^*$-algebras associated with hyperbolic dynamical systems were
studied
in~\cite{putnam,put_kam_sp,putnam:structure,ruelle:therm,ruelle:grpd}.

\section{Schur complements and rational multidimensional dynamics}
\subsection{Schur complements of operators}
Schur complements are widely used in linear algebra (usually
without knowledge of the name). The term ``Schur complement'' was
apparently introduced for the first time by
E.~V.~Haynsworth~\cite{haynsworth:inertia}. See a
survey~\cite{cottle:schursurvey} of the use of Schur complement in
algebra and statistics. It is also related to \emph{Bruhat normal
form} of matrices over a skew field (see~\cite{cohn}), which is
used to define the \emph{Dieudonn\'e determinant}.

Schur complements are quite often used in different situations
where renormalization principle can be used. Here we describe one
such situation, which is related to computation of spectra of
Hecke type operators and discrete Laplace operators attached to
self-similar groups and their Schreier graphs. The material
written here summarizes the ideas and results of the
articles~\cite{bgr:spec,grineksu_en,gr_zu:lamp,gri:solvedunsolved,grisavchuksunic:img}.

Let $H$ be a Hilbert space decomposed into a direct sum $H=H_1\oplus H_2$ of two non-zero
subspaces. Let $M\in B(H)$ be a bounded operator and let
\[M=\left(\begin{array}{cc}A & B\\ C & D\end{array}\right)\]
be an operator matrix representing $M$ according to this decomposition.

\begin{definition}
(i) Assume that $D\in B(H_2)$ is invertible. Then the \emph{first Schur complement}
is the operator
\[S_1(M)=A-BD^{-1}C\]

(ii) Assume that $A\in B(H_1)$ is invertible. Then the \emph{second Schur complement}
is the operator
\[S_2(M)=D-CA^{-1}B.\]
\end{definition}

The reason why Schur complements are useful for the spectral
problem is the following well known statement.

\begin{theorem}
\label{th:schurcomplement} Suppose that $D$ is invertible. Then
$M$ is invertible if and only if $S_1(M)$ is invertible.
Similarly, if $A$ is invertible, then $M$ is invertible if and
only if $S_2(M)$ is invertible.

The inverse is computed then by the formula
\begin{equation}\label{eq:frobenius} M^{-1}=\left(\begin{array}{cc}S_1^{-1} & -S_1^{-1}BD^{-1}\\
-D^{-1}CS_1^{-1} &
D^{-1}CS_1^{-1}BD^{-1}+D^{-1}\end{array}\right),\end{equation}
where $S_1=S_1(M)$.
\end{theorem}

\begin{proof}
Consider \[L=\left(\begin{array}{rr}1 & 0\\ -D^{-1}C & D^{-1}\end{array}\right),\]
where $0$ and $1$ represent the zero and identity linear maps between the corresponding subspaces of
$H$.

Then
\[L^{-1}=\left(\begin{array}{rr}1 & 0\\ C & D\end{array}\right),\quad ML=
\left(\begin{array}{rr} A-BD^{-1}C & BD^{-1}\\ 0 & 1\end{array}\right).\]
A triangular operator matrix \[R=\left(\begin{array}{rr} x & y\\ 0 & 1\end{array}\right)
\]
represents a (right) invertible operator if and only if $x$ is
(right) invertible, hence invertibility of $ML$ is equivalent to
invertibility of $S_1(M)$ and the result follows. The second part
of the theorem is proved similarly.

Direct computations show that~\eqref{eq:frobenius} takes place.
\end{proof}

Note that if $A$, $D$ and $M$ are invertible, then
\[M^{-1}=\left(\begin{array}{cc}S_1^{-1} & B'\\
C' & S_2^{-1}\end{array}\right),\] where $S_1=S_1(M)$,
$S_2=S_2(M)$ and
\begin{gather*}
B'=-S_1^{-1}BD^{-1}=-A^{-1}BS_2^{-1},\\
C'=-D^{-1}CS_1^{-1}=-S_2^{-1}CA^{-1}.
\end{gather*}
The last two equalities follow from
\[S_1A^{-1}B=B-BD^{-1}CA^{-1}B=BD^{-1}S_2\] and
\[CA^{-1}S_1=C-CA^{-1}BD^{-1}C=S_2D^{-1}C.\]

Formula~\eqref{eq:frobenius} is called sometimes \emph{Frobenius
formula}, see, for instance~\cite{gantmacher}.

We see that taking Schur complement of $M$ is equivalent to
inverting the left top corner of the matrix $M^{-1}$.

The following corollary is known as \emph{Schur
formula}~\cite{schur:formula} and easily follows from the proof of
Theorem~\ref{th:schurcomplement}.

\begin{corollary}
Let $H$ be finite dimensional and suppose that the determinant $|D|$ is not equal to zero.
Then \[|M|=|S_1(M)|\cdot |D|.\]
\end{corollary}

There is nothing special in decomposition of $H$ into a direct sum
of \emph{two} subspaces. If $H=H_1\oplus H_2\oplus\cdots\oplus
H_d$ and
\begin{equation}\label{eq:m}
M=\left(\begin{array}{ccc} M_{11} & \ldots & M_{1d}\\ \vdots &
\ddots & \vdots \\ M_{d1} & \ldots & M_{dd}\end{array}\right),
\end{equation}
where $M_{ij}:H_j\arr H_i$, then we can write $H$ as $H=H_1\oplus
H_1^{\perp}$, where $H_1^\perp=H_2\oplus\cdots\oplus H_d$ and get
\[M=\left(\begin{array}{cc}M_{11} & B\\ C & D\end{array}\right),\]
where $B, C$ and $D$ are operators represented by $1\times (d-1),
(d-1)\times 1$ and $(d-1)\times (d-1)$-sized operator matrices
coming from~\eqref{eq:m}. Then $S_1(M)$ is defined as before under
the condition that $D$ is invertible.

Changing the order of the summands (putting $H_i$ on the first place, say using a cyclic permutation)
we define the $i$th Schur complement $S_i(M)$.

We will use these definitions in situations when $H_i$ for $i=1,
\ldots, m$ are isomorphic. There are two cases.

(I) Finite dimensional, when $H=\underbrace{H'\oplus\cdots\oplus
H'}_{\text{$d$ times}}$ for $\dim H'<\infty$.

We will apply the Schur complements to the sequence of operators
$M_n\in B(H_n)$, where $\dim H_n=d^n$, and
\[H_{n+1}=\underbrace{H_n\oplus\cdots\oplus H_n}_{\text{$d$ times}},\]
for all $n=0, 1, \ldots$ and $\dim H_0=1$.

(II) Infinite dimensional case. Suppose that we have an infinite
dimensional (separable) Hilbert space $H$ and a fixed
$d$-similarity
\[\psi:H\arr H^d.\]

Let \[\widetilde{S}_i=\psi^{-1}S_i\psi:B(H)\arr B(H),\quad
i=1,\ldots, d\] be partially defined transformations, where $S_i$,
as before, is the $i$th Schur complement.

In the sequel we will write just $S_i$ instead of
$\widetilde{S}_i$, as this will not lead to a confusion.

\begin{proposition}
\label{pr:schurcuntz} Let $\rho$ be a representation of the Cuntz
algebra $O_d$ associated with the $d$-similarity $\psi:H\arr H^d$
and let $T_i=\rho(a_i)$ be the images of the generators of the
Cuntz algebra.  Then the Schur complements are given by
\[S_i(M)=\left(T_i^*M^{-1}T_i\right)^{-1}.\]
\end{proposition}

\begin{proof}
Follows directly from the definition of $T_i$ and inversion
formula~\eqref{eq:frobenius} from
Theorem~\ref{th:schurcomplement}. Recall (see
Proposition~\ref{pr:cuntzrepr}) that $T_i$ are given by
\[T_i(\xi)=\psi^{-1}(0, \ldots, 0, \xi, 0, \ldots, 0),\] where
$\xi$ is on the $i$th coordinate.
\end{proof}

\begin{corollary}
We have the following formula for the composition of Schur
complements
\[S_{i_1}\circ\cdots\circ S_{i_k}(M)=\left((T_{i_k}\cdots
T_{i_1})^*\cdot M^{-1}\cdot(T_{i_k}\cdots T_{i_1})\right)^{-1}.\]
\end{corollary}

We see that a composition of $k$ Schur complements associated with
the $d$-similarity $\psi:H\arr H^d$ is a Schur complement
associated with the corresponding $d^k$-similarity of $H$ obtained
by iteration of $\psi$.

Consider the semigroup $\mathcal{S}=\langle S_1, \ldots,
S_m\rangle$ of partial transformations of $B(H)$ generated by the
transformations $S_i$. We will call it the \emph{Schur
renorm-semigroup}. In the examples that will follow, we will
restrict $\mathcal{S}$ (or its particular elements
$s\in\mathcal{S}$) onto some $\mathcal{S}$-invariant (or
$s$-invariant) subspaces $B'\subset B(H)$. In the examples that we
consider $B'$ will be finite dimensional.

We will see in Example~\ref{ex:resolvent} how nicely Schur
complement behaves on the space of resolvents.

The Schur transformations $S_i$ are homogeneous:
\[S_i(tM)=tS_i(M)\]
for all $t\in\mathbb{C}$. This allows to define the corresponding
transformations $\widehat S_i$ on the ``projective space''
$PB(H)=B(H)/\C^{\times}$ of $B(H)$ (here $\C^\times$ is the
multiplicative group of $\C$). Restricting to invariant
finite-dimensional subspaces $B'\subset B(H)$ (real or complex) we
get transformations on the corresponding projective spaces
$\mathbb{R}P^n$ or $\mathbb{C}P^n$, where $n+1=\dim B'$.

\subsection{Schur complements in self-similar
algebras}\label{ss:scssalg} If $B$ is a unital Banach algebra
together with a unital embedding $B\arr M_{d\times d}(B)$ (we call
such algebras \emph{self-similar}), then we can compute the Schur
complements $S_i(a)$ of the elements of $B$.

Let $(G, \alb)$ be a self-similar group, let
$\phi:G\arr\symm\wr_\alb G$ be the corresponding wreath recursion
and let $\cpg$ be the universal Cuntz-Pimsner algebra of $G$. Then
the Schur complements $S_i$ are defined on $\maxalg\subset\cpg$
and can be computed using Proposition~\ref{pr:schurcuntz} or
directly using the matrix recursion $\phi:\C[G]\arr M_{d\times
d}(\C[G])$.

If $\psi:H\arr H^d$ is a fixed $d$-similarity on $H$ (i.e., a
representation of the Cuntz algebra $O_d$, then, as it was
observed in Proposition~\ref{pr:cpimsner}, every self-similar
representation $\rho$ of $G$ extends to a representation
$\rho:\cpg\arr B(H)$. The corresponding Schur complements on the
closure $\A_\rho$ of $\C[G]$ are also computed by the formula in
Proposition~\ref{pr:schurcuntz}.

Let us consider, for instance, the space $L^2(\xo, \nu)$, the
natural unitary representation $\pi$ of a self-similar group $(G,
\alb)$ on it and the natural $d$-similarity (and the associated
representation of the Cuntz algebra) on $L^2(\xo, \nu)$, see
Example~\ref{ex:murepr}. Consider the respective self-similar
completion $\A_{mes}$ of $\C[G]$.

The partition $\xi_n$ of the boundary $\partial\xs=\xo$ into a
disjoint union of $d^n$ cylinder subsets $v\xo$ ($v\in\alb^n$),
corresponding to the vertices of the $n$th level of the tree is
$G$-invariant (since it is invariant under the action of the whole
automorphism group of the tree).

Let $H_n=\{X_{A_i^{(n)}},\; 1\le i\le d^n\}$ be the finite
dimensional subspace spanned by the characteristic functions of
the atoms $A_i^{(n)}$ of the partition $\xi_n$. Then $H_n$ is a
$\pi(G)$-invariant subspace and $\{H_n\}_{n=1}^\infty$ is a nested
sequence with canonical embeddings $H_n\hookrightarrow H_{n+1}$
and
\begin{equation}\label{eq:Hn} H=\overline{\bigcup_{n=1}^\infty
H_n}.\end{equation}

Note that
\[H_{n+1}=\bigoplus_{1\le i\le d}T_i(H_n),\]
where $T_i=\pi(a_i)$ are the generators of the Cuntz algebra
associated to the natural $d$-similarity of $L^2(\xo, \nu)$. The
space $T_i(H_n)$ is the linear span of the characteristic
functions of cylindrical sets of the level $n$ with the first
letter $i$.

Let $M=\sum_{g\in G}\lambda_g g\in\ell^1(G)$ be an arbitrary
element of the Banach group algebra and let $\pi(M)$ be its image
in $\A_{mes}$. Denote $M_n=\pi(M)|_{H_n}$. Relation~\eqref{eq:Hn}
and $G$-invariance of the spaces $H_n$ implies that
\begin{equation}\label{eq:spec}
\spec{\pi(M)}=\overline{\bigcup_{n=1}^\infty\spec{M_n}}.\end{equation}

The formula~\eqref{eq:spec} reduces the problem of computation of
the spectrum of $M$ to finite dimensional problems of finding
$\spec{M_n}$.

The latter problem requires computation of the $d^n\times
d^n$-size matrices $g^{(n)}$ given by the recursive relations:
\[g^{(0)}=1,\quad g^{(n+1)}=\left(A_{yx}\right)_{x,y\in\alb},\quad
A_{yx}=\left\{\begin{array}{ll}g|_x^{(n)} & \text{if $g(x)=y$,}\\
0 & \text{otherwise,}\end{array}\right.\] repeating the matrix
recursion~\eqref{eq:matrrecursion} and finding spectra of their
linear combinations. There is no general tool for solving this
problem. Nevertheless, we will show how the problem can be solved
using the Schur map in some particular cases.

We will start with some examples which have a complete solution of
the spectral problem, then show a few examples when the spectral
problem is reduced to a problem in Dynamical Systems of finding
invariant sets for a multidimensional rational mapping (it looks
that these invariant sets are indeed ``strange attractors'' for
the corresponding maps, at least it is so in many cases). Finally
we will show some examples for which the spectral problem is
extremely interesting and has links with other topics, but it
looks like the method of Schur complements doesn't work for these
examples.

We perform some computations in homogeneous coordinates of the
space of parameters, since the spectrum of a pencil (i.e., the set
of values of parameters corresponding to degenerate matrices) is
invariant under multiplication by a non-zero number and the Schur
complement transformations are homogeneous. But the obtained
transformations are written in non-homogeneous coordinates as
transformations in the Euclidean space.

In all examples that follow we use parametric family of elements
of a group algebra (usually the sum with coefficients-parameters
of generators and identity) or a matrix with entries of this type
and consider the corresponding operator in $L^2(\xo, \nu)$ given
by the representation $\pi$ (so that the operators belong to the
algebra $\mathcal{A}_{mes}$).

\begin{example} Let $\mathfrak{G}=\langle a, b, c, d\rangle$ be the group from
Example~\ref{sss:grigorchuk}. Since the generators are of order 2,
the operator $M=\pi(a+b+c+d)$ is a self-adjoint operator in
$B(L^2(\partial T, \nu))$, being also a self-adjoint element of
$\A_{mes}$ with the same spectrum. In what follows we omit $\pi$.
We have $a=\left(\begin{array}{cc} 0 & 1\\ 1 &
0\end{array}\right)$, $b=\left(\begin{array}{cc} a & 0\\ 0 &
c\end{array}\right)$, $c=\left(\begin{array}{cc} a & 0\\ 0 &
d\end{array}\right)$, $d=\left(\begin{array}{cc} 1 & 0\\ 0 &
b\end{array}\right)$, and
\[M=a+b+c+d=\left(\begin{array}{cc} 2a+1 & 1\\ 1 &
b+c+d\end{array}\right).\] Also $a_{n+1}=\left(\begin{array}{cc} 0 & 1\\
1 & 0\end{array}\right)$, $b_{n+1}=\left(\begin{array}{cc} a_n & 0\\
0 & c_n\end{array}\right)$, $c_{n+1}=\left(\begin{array}{cc} a_n & 0\\
0 & d_n\end{array}\right)$, $d_{n+1}=\left(\begin{array}{cc} 1 & 0\\
0 & b_n\end{array}\right)$ are recursive relations determining
matrices of generators viewed as operators in the spaces $H_{n+1}$
and $H_n$ respectively (0 and 1 represent here the zero operators
and the identity operator of the corresponding size). Instead of
computing spectrum of $M_n=a_n+b_n+c_n+d_n$ directly, one can try
to find spectrum of the whole pencil $M_n(x, y, z,
u)=xa_n+yb_n+zc_n+ud_n$ of matrices, where $x, y, z, u\in\C$.

If we find spectrum of each of the pencils $M_n(x, y, z, u)$, we
will solve the problem of finding spectrum of the pencil $M(x, y,
z, u)=xa+yb+zc+ud$ of infinite dimensional operators, since
spectrum $M(x, y, z, u)$ is the closure of the unions of spectra
of $M_n(x, y, z, u)$.

The 5-dimensional space of operators $M(x, y, z, u,
v)=xa+yb+zc+ud+v\cdot 1$ is invariant with respect to the Schur
complements which in this case are given by the following
formulae.
\[\widehat S_1:\left(\begin{array}{c} x\\ y\\ z\\ u\\
v\end{array} \right)\mapsto\left(\begin{array}{c} z+y\\
\frac{x^2(2yzv-u(y^2+z^2-u^2+v^2))}{(y+z+u+v)(y+z-u-v)(y-z+u-v)(-y+z+u-v)}\\
\frac{x^2(2zuv-y(-y^2+z^2+u^2+v^2))}{(y+z+u+v)(y+z-u-v)(y-z+u-v)(-y+z+u-v)}\\
\frac{x^2(2yuv-z(y^2-z^2+u^2+v^2))}{(y+z+u+v)(y+z-u-v)(y-z+u-v)(-y+z+u-v)}\\
u+v+\frac{x^2(2yzu-v(y^2+z^2+u^2-v^2))}{(y+z+u+v)(y+z-u-v)(y-z+u-v)(-y+z+u-v)}\end{array}\right)\]
and
\[
\widehat S_2:\left(\begin{array}{c} x\\ y\\ z\\ u\\
v\end{array}\right)\mapsto\left(\begin{array}{c}\frac{x^2(y+z)}{(u+v+y+z)(u+v-y-z)}\\
u\\ y\\ z\\
v-\frac{x^2(u+v)}{(u+v+y+z)(u+v-y-z)}\end{array}\right)
\]

The complement $\widehat S_2$ is not defined when
$(y+z)a+(u+v)\cdot 1$ is not invertible, i.e., if
\begin{equation}\label{eq:s2not1}
y+z+u+v = 0
\end{equation}
or
\begin{equation}\label{eq:s2not2}
y+z-u-v = 0,
\end{equation}
as $\spec{a}=\{\pm 1\}$.

The complement $\widehat S_1$ is not defined for the same
conditions~\eqref{eq:s2not1} and~\eqref{eq:s2not2} and also the
conditions obtained from them by application of the cyclic
permutation $(u, y, z)\mapsto (z, u, y)\mapsto (y, z, u)$, as
\begin{multline*}yc+zd+ub+v\cdot 1\\ =y\left(\begin{array}{cc}a & 0\\ 0 &
d\end{array}\right)+z\left(\begin{array}{cc}1 & 0\\ 0 &
b\end{array}\right)+u\left(\begin{array}{cc}a & 0\\ 0 &
c\end{array}\right)+v\left(\begin{array}{cc}1 & 0\\ 0 &
1\end{array}\right)\\
=\left(\begin{array}{cc}(y+u)a+(z+v)1 & 0\\ 0 &
yd+zb+uc+v1\end{array}\right).
\end{multline*}

Note that in both cases the complement $\widehat S_i$ is not
defined if and only if the denominator in the corresponding
formulae is equal to zero.

Note that the third iteration of the map $\widehat S_2$ fixes the
three middle coordinates, hence we actually get a family of maps
$\C^2\arr\C^2$ depending on three parameters.

It is not known at the moment if the spectrum of the pencil $M(x,
y, z, u, v)$ is invariant with respect to these maps and what are
its topological properties. Further investigation in this
direction would be interesting.

But let us simplify the problem and consider the pencil
$M(x)=xa+b+c+d$. The spectral problem for $M(x)$ consists in
finding those values of the parameters $x$ and $\lambda$ for which
the operator $M(x)-\lambda\cdot 1$ is not invertible. We will call
this set the \emph{spectrum} of the pencil. This terminology will
be used also for further examples. To simplify transformations let
us change the system of coordinates and consider the pencil
\[R(\lambda, \mu)=-\lambda a+b+c+d-(\mu+1)1\]
and its finite dimensional approximations
\[R_n(\lambda, \mu)=-\lambda a_n+b_n+c_n+d_n-(\mu+1)1_n\]
represented by a 2-parametric family of $2^n\times 2^n$-matrices.

Let \[\Sigma=\{(\lambda, \mu)\;:\;\text{$R(\lambda, \mu)$ is not
invertible}\}\] and
\[\Sigma_n=\{(\lambda, \mu)\;:\;\text{$R_n(\lambda, \mu)$ is not
invertible}\}.\] Define $F, G:\R^2\arr\R^2$ be two rational maps
\begin{eqnarray*}
F:\left(\begin{array}{c}\lambda\\
\mu\end{array}\right) &\mapsto&
\left(\begin{array}{c}\frac{2(4-\mu^2)}{\lambda^2}\\
-\mu-\frac{\mu(4-\mu^2)}{\lambda^2}\end{array}\right),\\
G:\left(\begin{array}{c}\lambda\\
\mu\end{array}\right) &\mapsto&
\left(\begin{array}{c}\frac{2\lambda^2}{4-\mu^2}\\
\mu+\frac{\mu\lambda^2}{4-\mu^2}\end{array}\right) \end{eqnarray*}
and observe, following Y.~Vorobets, that $H\circ F=G$, $H\circ
G=F$, where
\[H:\left(\begin{array}{c}\lambda\\
\mu\end{array}\right)\mapsto\left(\begin{array}{c}4/\lambda\\
-2\mu/\lambda\end{array}\right)\] is an involutive map
$\R^2\arr\R^2$.\end{example}

\begin{figure}
 \includegraphics{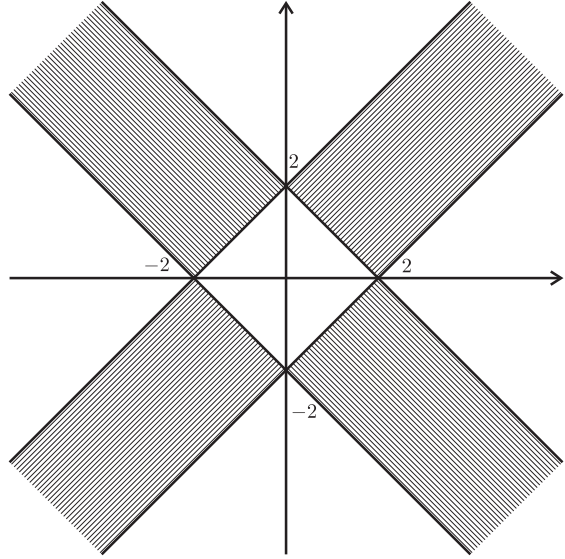}\\
 \caption{The spectrum $\Sigma$}\label{fig:sigma}
\end{figure}

\begin{figure}
\includegraphics{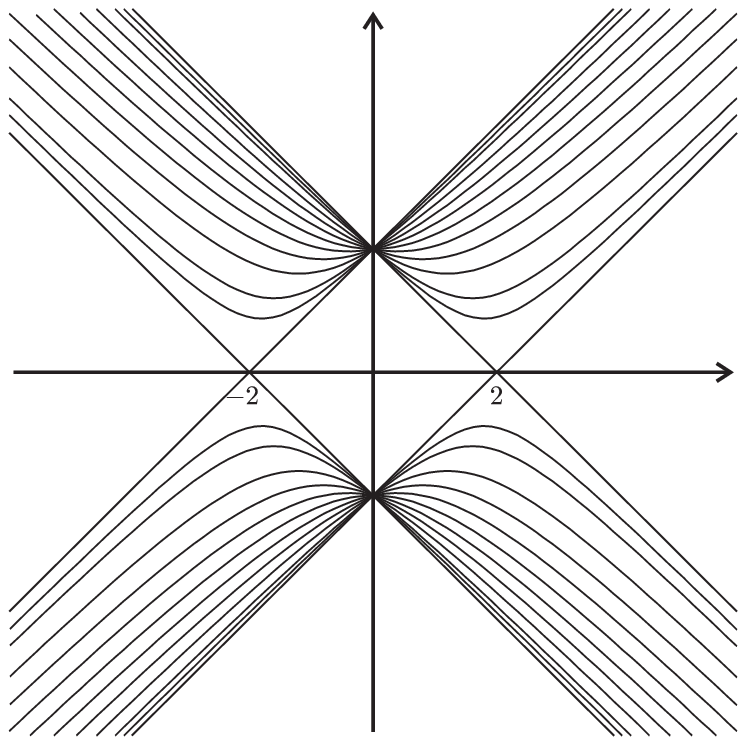}\label{fig:sigma5}
\caption{The spectrum $\Sigma_5$}
\end{figure}

\begin{theorem}
I. The maps $F$ and $G$ seen as maps on the projective space are
respectively the first and the second Schur maps restricted on the
pencil $R(\lambda, \mu)$.

II. The set $\Sigma$ is invariant with respect to $F$ and $G$,
i.e., $F^{-1}\Sigma=\Sigma$, $G^{-1}\Sigma=\Sigma$ and relations
\[\Sigma_{n+1}=F^{-1}\Sigma_n=G^{-1}\Sigma_n\]
hold.

III. The set $\Sigma$ is the set shown on Figure~\ref{fig:sigma}.

IV. The set $\Sigma_n$ is the union of the line $\lambda+\mu-2=0$
and hyperbolas $H_\theta=0$,
$\theta\in\bigcup_{i=0}\alpha^{-i}(0)$, where
$H_\theta=4-\mu^2+\lambda^2+4\lambda\theta$. See
Figure~\ref{fig:sigma5} for the spectrum $\Sigma_5$.
\end{theorem}

\begin{proof}
We restrict ourselves here only to computation of the Schur maps.
The rest can be found in~\cite{bgr:spec}.

We have
\[R(\lambda, \mu)=\left(\begin{array}{cc}2a-\mu & -\lambda \\ -\lambda &
b+c+d-\mu-1\end{array}\right).\]

For getting the first Schur map observe that $t=(b+c+d-1)/2$ is an
idempotent. Hence
\[(b+c+d-\mu-1)^{-1}=(2t-\mu)^{-1}=\frac{2t+\mu}{4-\mu^2}\]
and \[\widehat S_1(R(\lambda,
\mu))=2a-\mu-\frac{\lambda^2(b+c+d+\mu-1)}{4-\mu^2},\] which is
proportional to  $-\lambda'a+b+c+d-(\mu'+1)$ for $(\lambda',
\mu')=F(\lambda, \mu)$.

Similarly
\begin{multline*}
\widehat S_2(R(\lambda, \mu))=b+c+d-\mu-1-\lambda^2(2a-\mu)^{-1}\\
=b+c+d-\mu-1-\frac{\lambda^2(2a+\mu)}{4-\mu^2}\\
=-\frac{2\lambda^2}{4-\mu^2}a+b+c+d-\left(\mu+\frac{\mu\lambda^2}{4-\mu^2}+1\right)\\
=-\lambda''a+b+c+d-(\mu''+1),
\end{multline*}
where $(\lambda'', \mu'')=G(\lambda, \mu)$.
\end{proof}

\begin{remark}
The equations of hyperbolas $H_\theta$ are related by the rule
\[H_\theta(F(\lambda, \mu))=\frac{1-\mu^2}{\lambda^4}H_{\frac
12\sqrt{2-2\theta}}(\lambda, \mu) H_{-\frac
12\sqrt{2-2\theta}}(\lambda, \mu).\]

The maps $F$ and $G$ are semi-conjugate to the map $x\mapsto
2x^2-1$ via the maps $\psi_F(x, y)=\frac{4-y^2+x^2}{4x}$ and
$\psi_G(x, y)=\frac{4-x^2+y^2}{4y}$, respectively. This is the
crucial point in computation of the spectrum of $R_n(\lambda,
\mu)$.
\end{remark}

\begin{example}
Let $H=H^{(3)}$ be the Hanoi Towers group (on three pegs) studied
in~\cite{grisunik:hanoi} and~\cite{grisunik:hanoioberwolfach}. It
is a self-similar group acting on $\xs$ for $\alb=\{0, 1, 2\}$
with generators $a, b, c$ satisfying the following matrix
recursions
\begin{align*} a &=\left(\begin{array}{ccc}0 & 1 & 0\\ 1 & 0 & 0\\ 0 & 0 &
a\end{array}\right),\\
b &=\left(\begin{array}{ccc}0 & 0 & 1\\ 0 & b & 0\\ 1 & 0 &
0\end{array}\right),\\
c &=\left(\begin{array}{ccc}c & 0 & 0\\ 0 & 0 & 1\\ 0 & 1 &
0\end{array}\right).
\end{align*}
 Consider the two-parametric family of elements of the
 $C^*$-algebra $\mathcal{A}_{mes}$.
\begin{multline*}
\Delta(x, y)=\left(\begin{array}{ccc}c-x & y & y\\
y & b-x & y\\ y & y & a-x\end{array}\right)=\\
\left(\begin{array}{ccc|ccc|ccc}c-x & 0 & 0 & y & 0 & 0 & y & 0 &
0\\
0 & -x & 1 & 0 & y & 0 & 0 & y &
0\\
0 & 1 & -x & 0 & 0 & y & 0 & 0 & y\\ \hline
y & 0 & 0 & -x & 0 & 1 & y & 0 & 0\\
0 & y & 0 & 0 & b-x & 0 & 0 & y & 0\\
0 & 0 & y & 1 & 0 & -x & 0 & 0 & y\\ \hline
y & 0 & 0 & y & 0 & 0 & -x & 1 & 0\\
0 & y & 0 & 0 & y & 0 & 1 & -x & 0\\
0 & 0 & y & 0 & 0 & y & 0 & 0 & a-x
\end{array}\right).
\end{multline*}
Permuting rows and columns and dividing them into blocks we get
the matrix
\[\left(\begin{array}{ccc|cccccc}
c-x & 0  & 0  & y & 0 & 0 & y & 0 & 0 \\
0  & b-x & 0  & 0 & y & 0 & 0 & y & 0 \\
0  & 0  & a-x & 0 & 0 & y & 0 & 0 & y \\ \hline
y  & 0  & 0  & -x & 0 & 1 & y & 0 & 0 \\
0  & y  & 0  & 0 & -x & 0 & 0 & y & 1 \\
0  & 0  & y  & 1 & 0 & -x & 0 & 0 & y \\
y  & 0  & 0  & y & 0 & 0 & -x & 1 & 0 \\
0  & y  & 0  & 0 & y & 0 & 1 & -x & 0 \\
0  & 0  & y  & 0 & 1 & y & 0 & 0 & -x
\end{array}\right).\]
Computation of Schur complement with respect to the given
partition of the matrix yields
\[S_1(\Delta(x, y))=\Delta(x', y'),\]
where
\[x'=x-\frac{2(x^2-x-y^2)y^2}{(x-y-1)(x^2-1+y-y^2)}\]
and
\[y'=\frac{(x+y-1)y^2}{(x-y-1)(x^2-1+y-y^2)}.\]
These rational functions were calculated in~\cite{grisunik:hanoi}
using a different base. As it was observed
in~\cite{grisunik:hanoi}, the map $F:(x, y)\mapsto (x', y')$ is
semiconjugate to the map $f:\R\arr\R:x\mapsto x^2-x-3$. The
spectrum of $\Delta(x, y)$ is in this case the union
$\bigcup_{\theta\in\bigcup f^{-n}(S)}\mathcal{H}_\theta\cup
L_0\cup L_1\cup L_2$, where $S=\{-2, 0\}$, and
$\mathcal{H}_\theta$ is the hyperbola $x^2-xy-2y^2-\theta y=1$ and
$L_0, L_1, L_2$ are the lines given by the equations
\begin{align*}
x-1-2y &= 0,\\
x+1+y &= 0,\\
x-1+y &= 0.
\end{align*}
A part of the spectrum is drawn on Figure~\ref{fig:hanoi}.

\begin{figure}[h]
\centering\includegraphics{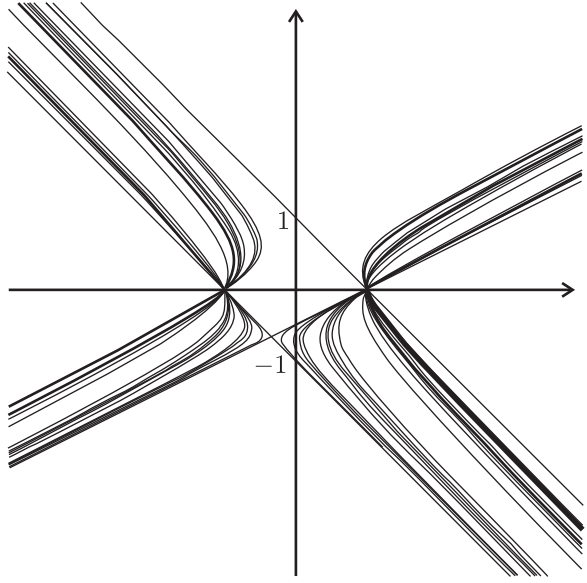} \caption{Spectrum of
$H_3$}\label{fig:hanoi}
\end{figure}

\end{example}

\begin{example}
Let $B=\img{z^2-1}$ be the Basilica group, studied
in~\cite{zukgrigorchuk:3st,zukgrigorchuk:3stsp,barthvirag,bgn,nek:book}.
It is realized as a self-similar group acting on the binary tree
and generated by $a, b$, which are given by the matrix recursions
\[a=\left(\begin{array}{cc}1 & 0\\ 0 & b\end{array}\right),\quad
b=\left(\begin{array}{cc}0 & 1\\ a & 0\end{array}\right).\]
Consider the pencil
\[R(\lambda, \mu)=a+a^{-1}+\lambda(b+b^{-1})-\mu=
\left(\begin{array}{cc}2-\mu & \lambda(1+a^{-1})\\
\lambda(1+a) & b+b^{-1}-\mu\end{array}\right).\]

We have
\begin{multline*}
\widehat S_2(R(\lambda,
\mu))=b+b^{-1}-\mu-\frac{\lambda^2(2+a+a^{-1})}{2-\mu}\\
=\frac{\lambda^2}{\mu-2}(a+a^{-1})+b+b^{-1}-\mu+\frac{2\lambda^2}{\mu-2}\\
\sim
a+a^{-1}+\frac{\mu-2}{\lambda^2}(b+b^{-1})-\left(\frac{\mu(\mu-2)}{\lambda^2}-2\right),
\end{multline*}
where ``$\sim$'' means ``proportional''. We get hence the rational
map
\begin{equation}\label{eq:basilica}
\left\{\begin{array}{rcl}\lambda & \mapsto &
\frac{\mu-2}{\lambda^2},\\
\mu &\mapsto & -2+\frac{\mu(\mu-2)}{\lambda^2}.\end{array}\right.
\end{equation}

This rational map is quite complicated and the structure of the
spectrum of the pencil $R(\lambda, \mu)$ is unclear. Computer
experiments suggest that it has to have a structure of a ``strange
attractor''. Most probably the map~\eqref{eq:basilica} is not
semiconjugate to any one-dimensional map.

The pencil $R(\lambda, \mu)$ is not invariant with respect to the
first Schur complement (because the inverse of $b+b^{-1}-\mu$ is
not a finite sum).
\end{example}

\begin{example}
Let $G=\img{z^2+i}$ be the group studied
in~\cite{bgn,nek:book,grisavchuksunic:img}. It is generated by $a,
b, c$ given by the matrix recursions
\[a=\left(\begin{array}{cc}0 & 1\\ 1 & 0\end{array}\right),\quad
b=\left(\begin{array}{cc}a & 0\\ 0& c\end{array}\right),\quad
c=\left(\begin{array}{cc}b & 0\\ 0 & 1\end{array}\right).\] It is
a branch group of intermediate growth. For the notion of branch
groups see~\cite{handbook:branch} and~\cite{grigorchuk:branch};
the proof of intermediate growth of $\img{z^2+i}$ is given
in~\cite{buxperez:imgi}. Consider the pencil
\[M(y, z, \lambda)=a+yb+zc-\lambda=M(y, z, \lambda)=\left(\begin{array}{cc}
ya+zb-\lambda & 1\\ 1 & yc+z-\lambda\end{array}\right),\] then
\[\widehat S_1(M(y, z, \lambda)=M(y', z', \lambda'),\] where $\Phi:(y, z,
\lambda)\mapsto (y', z', \lambda')$ is the map
\[\Phi:\left\{\begin{array}{rcl} y & \mapsto & \frac{z}{y},\\
z & \mapsto &
\frac{1}{(z-\lambda-y)(z-\lambda+y)},\\
\lambda & \mapsto & \frac{-\lambda
y^2+\lambda(z-\lambda)^2+z-\lambda}{y(z-\lambda-y)(z-\lambda+y)}.\end{array}\right.\]
It is not known if $\Phi$ is semiconjugate to a map of a smaller
dimension and therefore there is no information about the
structure of the spectrum of the pencil $M(y, z, \lambda)$.

In this example, as also in the previous one, we have difficulty
to apply $S_2$, since in these cases we do not get finite
combinations of the group elements.
\end{example}

\begin{example}
Consider now the free group from Example~\ref{sss:free}. The
corresponding matrix recursion is
\[a_{n+1}=\left(\begin{array}{cc} b_n & 0\\ 0 &
b_n\end{array}\right),\quad b_{n+1}=\left(\begin{array}{cc} 0 & c_n\\
a_n & 0\end{array}\right),\quad c_{n+1}=\left(\begin{array}{cc} 0 & a_n\\
c_n & 0\end{array}\right).\] For the inverses we have
\[a_{n+1}^{-1}=\left(\begin{array}{cc} b_n^{-1} & 0\\ 0 &
b_n^{-1}\end{array}\right),\quad b^{-1}_{n+1}=\left(\begin{array}{cc} 0 & a_n^{-1}\\
c_n^{-1} & 0\end{array}\right),\quad c^{-1}_{n+1}=\left(\begin{array}{cc} 0 & c_n^{-1}\\
a_n^{-1} & 0\end{array}\right).\] The problem is to find the
spectrum of the matrix
\begin{multline*}M_n=a_n+a^{-1}_n+b_n+b_n^{-1}+c_n+c_n^{-1}=\\ \left(\begin{array}{cc}
b_{n-1}+
b^{-1}_{n-1} & a_{n-1}+a^{-1}_{n-1}+c_{n-1}+c^{-1}_{n-1}\\
a_{n-1}+a^{-1}_{n-1}+c_{n-1}+c^{-1}_{n-1} &
b_{n-1}+b^{-1}_{n-1}\end{array}\right),\end{multline*} but
introduction of new parameters and computation of Schur
complements unfortunately does not lead us to a success.
\end{example}

\begin{question} What is the spectrum of the matrix $M_n$? Does there exist
$\epsilon>0$ such that for every $n$ the spectrum of $M_n$ belongs
to the set $[-6+\epsilon, 6-\epsilon]\cup\{6, -6\}$? If the answer
is ``yes'', then the sequence of the Schreier graphs of the action
of $\langle a, b, c\rangle$ on the levels of the tree $\xs$ is a
sequence of expanders. (Here a \emph{Schreier graph} of a group
generated by a finite set $S$ and acting on a set $M$ is the graph
with the set of vertices $M$ in which two vertices $x, y\in M$ are
connected by an edge if one is the image of the other under the
action of an element of $S$.)
\end{question}

\begin{example}
Consider the group generated by the transformations $a, b, c$ of
$\xs$ for $\alb=\{0, 1\}$ given by the recurrent relations
\begin{alignat*}{2}
a(0w) &= 1b(w), &\quad a(1w) &= 0b(w),\\
b(0w) &= 0a(w), &\quad b(1w) &= 1c(w),\\
c(0w) &= 0c(w), &\quad c(1w) &= 1a(w),
\end{alignat*}
i.e., the group generated by the automaton shown on
Figure~\ref{fig:freepr}.

\begin{figure}
 \includegraphics{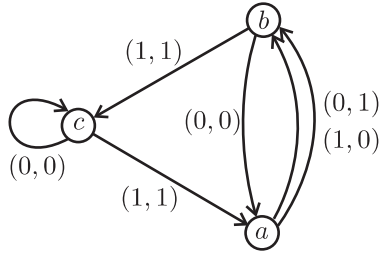}\\
 \caption{An automaton generating $C_2*C_2*C_2$}\label{fig:freepr}
\end{figure}

This group is isomorphic (by a result of Y.~Muntyan and
D.~Savchuk) to the free product $C_2*C_2*C_2$ of three groups of
order 2, see~\cite{nek:book} Theorem~1.10.2. The elements $a, b,
c$ are of order 2.

The corresponding matrix recursions are \[
a_{n+1}=\left(\begin{array}{cc} 0 & b_n \\ b_n &
0\end{array}\right),\quad b_{n+1}=\left(\begin{array}{cc} a_n & 0
\\ 0 & c_n\end{array}\right),\quad c_{n+1}=\left(\begin{array}{cc} c_n & 0
\\ 0 & a_n\end{array}\right).\]
\end{example}

\begin{question}
What is the spectrum of the matrices
\[M_n=a_n+b_n+c_n=\left(\begin{array}{cc} a_{n-1}+c_{n-1} & b_n \\ b_n &
a_{n-1}+c_{n-1}\end{array}\right)?\] In particular, is the
spectrum of $M$ a subset of $(-2\sqrt{2}, 2\sqrt{2})\cup\{3\}$? If
it is, then the Schreier graphs of the action of the group on the
levels of the tree $\xs$ are Ramanujan.
\end{question}

\section{Analytic families}

Let $\psi:H\arr H_1\oplus H_2$ be an isomorphism and let
\[M(z)=\left(\begin{array}{cc}A(z) & B(z)\\ C(z) &
D(z)\end{array}\right)\] be an analytic on some domain $\Omega$
operator valued function. Assume that $D(z)$ s invertible on
$\Omega$. Then \[S_1(M(z))=A(z)-B(z)D^{-1}(z)C(z)\] is analytic on
$\Omega$. We get therefore a partially defined Schur map between
the spaces of analytic operator valued functions.

If $H$ is a Hilbert space with a $d$-similarity $\psi:H\arr H^d$,
then we can define in a similar way partially defined Schur
transformations $\wt S_i$ on the space of analytic operator valued
functions on $H$.

\begin{example}\label{ex:resolvent} Let $M=\left(\begin{array}{cc}A & B\\ C &
D\end{array}\right)$ and consider the function
\[M-z=\left(\begin{array}{cc}A-z & B\\ C & D-z\end{array}\right).\]
It is mapped by the Schur map to $S_1=A-z-B(D-z)^{-1}C$ and we
have
\[(M-z)^{-1}=\left(\begin{array}{cc}S_1^{-1} & -S_1^{-1}B(D-z)^{-1}\\
-(D-z)^{-1}CS_1^{-1} &
-(D-z)^{-1}CS_1^{-1}B(D-z)^{-1}+(D-z)^{-1}\end{array}\right).\] It
follows from~\eqref{eq:frobenius} in
Theorem~\ref{th:schurcomplement} that $S_1((M-z)^{-1})=(A-z)^{-1}$
and, similarly, $S_2((M-z)^{-1})=(D-z)^{-1}$.
\end{example}

Consider now an operator-valued holomorphic function of $n$
complex variables
\[M:\C^n\arr B(H)\]
and let again $\psi:H\arr H^d$ be a $d$-similarity. We will denote
the respective Schur transformations $\wt S_i$ just by $S_i$.

\begin{definition}
The function $M(z)$ is \emph{self-similar} (with respect to $\wt
S_i$) if there is a function $F:\C^n\arr\C^n$ such that
\[S_i(M(z))=M(F(z)).\]
\end{definition}

We can consider functions $M$ from the projective space $\C
P^{n-1}$ to $B(H)$. Self-similarity of such functions is defined
analogically.

\begin{example}
Let $G$ be a self-similar group generated by $\{a_1, \ldots,
a_k\}$ and acting on a rooted $d$-regular tree $\xs$. Consider the
pencil of Hecke type operators
\[M(z)=M(z_1, \ldots, z_k)=\sum_{i=1}^kz_i\pi(a_i+a_i^{-1}),\]
where $\pi$ is the self-similar representation of $G$ on $L^2(\xo,
\nu)$. The operators $M(z)$ belong to the algebra $\A_{mes}$.

Assume that $M(z)$ is a self-similar function with respect to
$S_i$. This means that
\[S_i(M(z_1, \ldots, z_k))=\sum_{i=1}^kz_i'\pi(a_i+a_i^{-1})\]
where $z_j'=Z_j(z_1, \ldots, z_k)$ for $j=1, \ldots, k$ are some
functions. We get hence a map $F=(Z_1, \ldots, Z_k):\C^k\arr\C^k$
and the projectivized map $\wt F:\C P^{k-1}\arr\C P^{k-1}$.
\end{example}

\begin{proposition}
\label{pr:Linvariant} The hyperspace
$L=\left\{\sum_{i=1}^kz_i=0\right\}$ is forward $S_i$-invariant
for all $i$.
\end{proposition}

Here forward $S_i$-invariance is the condition
$S_i(L\cap\mathop{\mathrm{Dom}}(S_i))\subset L$.

\begin{proof}
\begin{lemma}
If $\left(\begin{array}{c}x_1\\ x_2\end{array}\right)$ belongs to
the kernel of $M$, then $x_1\in\ker S_1(M)$. Conversly, if
$x_1\in\ker S_1(M)$ then there is $x_2$ such that
$\left(\begin{array}{c}x_1\\ x_2\end{array}\right)\in\ker M$.
\end{lemma}

\begin{proof}
We have \[\left(\begin{array}{cc}A & B \\ C &
D\end{array}\right)\left(\begin{array}{c} x_1\\
x_2\end{array}\right)=\left(\begin{array}{c}Ax_1+Bx_2 \\
Cx_1+Dx_2\end{array}\right)=\left(\begin{array}{c}0\\
0\end{array}\right).\] We get $x_2=-D^{-1}Cx_1$ if $D$ has right
inverse. Consequently, $(A-BD^{-1}C)x_1=0$.

Conversly, if $x_1\in\ker S_1(M)$, then
$\left(\begin{array}{c}x_1\\ x_2\end{array}\right)$, where
$x_2=-D^{-1}Cx_1$.
\end{proof}

Since $\sum_{i=1}^kz_i=0$ the constant function $c\ne 0$ on
$\partial T$ belongs to $\ker M$. Its restriction onto the
cylindrical set of words starting with $i$ (where $i$ is the same
as in $S_i$) is also a constant function and by the above lemma,
we have $S_i(M)c=0$, hence $\sum z_i'=0$. That means that
$S_i(L\cap\mathop{\mathrm{Dom}}(S_i))\subset L$.
\end{proof}

Let \[\Sigma M(z)=\{z\;:\;M(z)\text{\ is not invertible}\}\] be
the spectrum (critical set) of the pencil $M(z)$. We have then the
following corollary of Theorem~\ref{th:schurcomplement}

\begin{corollary}
Let $M(z)$ and $A(z), B(z), C(z), D(z)$ be as before. Then
\[\Sigma M(z)\backslash\Sigma D(z)=\Sigma S_1(M(z))\setminus \Sigma D(z).\]
\end{corollary}

If $M(z)$ is self-similar, i.e., if $S_1(M(z))=M(F(z))$, then
\[\Sigma S_1(M(z))=\Sigma M(F(z))=F^{-1}(\Sigma(M(z))),\]
and we get

\begin{corollary}
The spectrum $\Sigma M(z)$ is backward-invariant under $F$, i.e.,
\[F^{-1}(\Sigma M(z))=\Sigma M(z)\]
if and only if \[\Sigma D(z)\cap \Sigma M(z)=\Sigma D(z)\cap
F^{-1}(\Sigma M(z)).\]
\end{corollary}

The condition of the corollary is not easy to check if
$\Sigma=\Sigma M(z)$ is unknown. In all examples that were
treated, though, we have $F^{-1}\Sigma=\Sigma$.

\begin{question} Under what natural and easy-to-check
conditions the equality $F^{-1}\Sigma=\Sigma$ is true?
Under what conditions we have
$\Sigma=\overline{\bigcup_{n=0}^\infty F^{-n}L}$?
\end{question}

In all treated examples the map $F$ is rational. Therefore the
problem of finding the spectrum of a pencil in a self-similar
group is related to the problem of description of invariant
subsets of multidimensional rational mappings. This subject is of
independent interest in dynamical systems, multidimensional
complex analysis, etc (see~\cite{sibony}). Spectra of Hecke type
elements in $C^*$-algebras related to self-similar groups is a big
source of interesting examples of dynamical systems on the complex
projective space.

\section{Schur maps and random walks}
Consider the map $J:B(H)\arr B(H):A\mapsto A+I$ where $I$ is the
identity operator and suppose that we have fixed a $d$-similarity
$\psi:H\arr H^d$. Consider the conjugates $k_i=JS_iJ^{-1}$ of the
Schur maps. We call $k_i$ the \emph{probabilistic Schur maps} (the
reason will be clarified later). If $M=\left(\begin{array}{cc}A &
B\\ C & D\end{array}\right)$, then
\begin{equation}\label{eq:k_1}k_1(M)=A+B(I-D)^{-1}C,\end{equation}
where $I$ is the identity operator (or a matrix of the same size
as $D$).

Following Bartholdi, Virag~\cite{barthvirag} and
Kaimanovich~\cite{kaimanovich:munchhausen}, we are going to apply
the probabilistic Schur maps to study random walks on self-similar
groups.

Let $(G, \alb)$ be a self-similar group acting on a $d$-regular
tree and let $\mu$ be a probability distribution on $G$. Thus,
every element $g$ of $G$ has a mass $\mu_g$, $0\le \mu_g\le 1$ and
$\sum_{g\in G}\mu_g=1$. The set $\supp\mu=\{g\;:\;\mu_g\ne 0\}$ is
called the \emph{support} of $\mu$. The measure $\mu$ is
\emph{non-degenerate} if $\supp\mu$ generates $G$. We will
identify $\mu$ with the element
\[\mu=\sum_{g\in G}\mu_g g\] of $\ell^1(G)$. Moreover, $\mu$
belongs to the simplex
\[\ell_+^1(G)=\{\mu\in\ell^1(G)\;:\;0\le\mu_g\le 1, \sum_{g\in G}\mu_g\}.\]
The left random walk on $G$ generated by $\mu$ starts at some
element of $G$ (usually at the identity) and at each step makes
the move $g\mapsto hg$ with probability $\mu(h)$.

When started at $e$ (the identity element of $G$) the distribution
at the moment $n$ is given by the $n$th convolution
$\mu^{(n)}=\underbrace{\mu*\cdots *\mu}_n$ of $\mu$ and the
probability of return $P_{e, e}^{(n)}$ is equal to $\mu^{(n)}(e)$.
Left random walk on $G$ is invariant with respect to the right
action of $G$ on itself.

The main topics that interest specialist in random walks are:
\begin{enumerate}
\item norm of the Markov operator \[M=\sum_{g\in G}\mu_gL_g\]
in $\ell^2(G)$, where $L_g$ are the operators of the left regular
representation. We identify here the elements of $\ell^1(G)$ with
the left convolution operators on $\ell^2(G)$;

\item spectrum of $M$;

\item asymptotic behavior of $P_{e, e}^{(n)}$ when $n\to \infty$;

\item Liouville property (i.e., when all harmonic functions are
constant);
\end{enumerate}
and other asymptotic characteristics of the random walks.

Let $u$ be one of the vertices of the first level of the tree and
let $H=\st_G(u)$ be its stabilizer in $G$. As $H$ has finite index
in $G$, the random walk hits $H$ with probability 1. Let $\mu_H$
be the distribution on $H$ given by the probability of the first
hit, i.e.,
\[\mu_H(h)=\sum_{n=0}^\infty f_{e, h}^{(n)},\]
where $f_{e, h}^{n}$ is the probability to hit $H$ at the element
$h$ for the first time at step $n$. As $H<G$ is a recurrent set as
a subgroup of finite index in $G$, we have $\sum_{h\in
H}\mu_H(h)=1$. Let $p_i:H\arr G$ be the $i$th projection map
$h\mapsto h|_i$ for $1\le i\le d$, and let $\mu_i$ be the image of
$\mu_H$ under $p_i$.

The next theorem and its proof are analogous to Theorem 2.3 of
V.~Kaimanovich~\cite{kaimanovich:munchhausen}, but they are
formulated a bit differently.

\begin{theorem}
In the above conditions
\[\mu_i=k_i(\mu).\]
\end{theorem}

\begin{proof}
Let $M=(m_{ij})_{1\le i, j\le d}$ be the matrix representation of
the Markov operator $M$ coming from the wreath recursion
$\phi:G\hookrightarrow \symm\wr_\alb G$. Let us extend it to the
map $\ell_+^1(G)\arr M_d(\ell_+^1(G))$.

Then the measures $\mu_i$ can be expressed as
\begin{equation}\label{eq:mii}\mu_i=m_{ii}+M_{i\overline i}\left(I+M_{\overline i\overline
i}+ M_{\overline i\overline i}^2+\cdots\right)M_{\overline i
i}=\mu_{ii}+M_{i\overline i}\left(I-M_{\overline i\overline
i}\right)^{-1}M_{\overline ii},\end{equation} where $M_{i\overline
i}$ (respectively, $M_{\overline i i}$) denotes the row
$(m_{ij})_{j\ne i}$ (respectively, the column $(m_{ji})_{j\ne i}$)
of the matrix $M$ with deleted element $m_{ii}$, and $M_{\overline
i\overline i}$ is the $(d-1)\times (d-1)$-matrix obtained from $M$
by removing its $i$th row and $i$th column.

The first term in~\eqref{eq:mii} corresponds to staying at the
point $i$ (in the random walk on $X$ induced by the random walk on
$G$), while the first factor of the second term corresponds to
moving from $i$ to $\alb\setminus\{i\}$, the second its factor
corresponds to staying in $\alb\setminus\{i\}$ and the third
factor of the second term corresponds to moving back from
$\alb\setminus\{i\}$ to $i$.

Comparing~\eqref{eq:mii} with~\eqref{eq:k_1} we get the statement
of the theorem.
\end{proof}

In case $\st_G(i)=\st_G(1)$ for all $i\in\alb$ (for instance, if
$\xs$ is a binary tree, or if $G$ acts on the first level by
powers of a transitive cycle $\sigma$) we can interpret the above
fact as that we have a sequence of stopping times $\tau(n)$ such
that if
\[\phi(Y^{(n)})=\sigma^{(n)}\left(Y^{(n)}_1, Y^{(n)}_2,\ldots,
Y^{(n)}_d\right)\] is image of the $n$th step $Y^{(n)}$ of the
random walk under the wreath recursion $\phi$, then
\[\phi(Y^{\tau(n)})=\left(Y^{(\tau(n))}_1, Y^{(\tau(n))}_2, \ldots Y^{(\tau(n))}_d\right),\]
i.e., the random element at time $\tau(n)$ belongs to the
stabilizer of the first level, and $Z_i^{(n)}=Y_i^{(\tau(n))}$ is
the random walk on $G$ determined by the measure $\mu_i$. Thus we
can treat asymptotic characteristics of $\mu$-random walk on $G$
via $\mu_i$-random walks on $G$. Of course for complete
reconstruction of $(G, \mu)$ we also need to know the joint
distributions of the processes $(G, \mu_i)$ (as they are usually
not independent). Nevertheless, some information about the random
walk can be obtained without the independence.

The maps $k_i:\mathcal{M}\arr\mathcal{M}:\mu\mapsto\mu_i$ on the
simplex of measures on $G$ are continuous and their fixed points
are of a special interest for us. The fixed points always exist
(for example the unit mass concentrated on the identity), but we
are interested only in non-degenerate fixed points (i.e., with
support generating the group).

\begin{remark}
Invariance of the hyperplane $L=\{x\in\ell^1(G)\;:\;\sum_{g\in
G}x_g=0\}$ under the Schur map (the proof of which is analogous to
the proof of Proposition~\ref{pr:Linvariant}) implies
$k_i$-invariance of the hyperplane $L=\{\lambda\;:\;\sum_{g\in
G}\lambda_g=1\}$. The probabilistic meaning of the maps $k_i$
shows that the simplex of probability measures $\mathcal{M}\subset
L$ is also $k_i$-invariant.
\end{remark}

Now we are going to describe briefly V.~Kaimanovich's approach to
testing amenability of $G$ by entropy method (which develops the
ideas previously expressed in~\cite{barthvirag}).

Having a left random walk $g_{n+1}=h_{n+1}g_n$, given by $(G,
\mu)$, where $(h_n)_{n\ge 0}$ is a sequence of independent
$\mu$-distributed random variables, we consider the induced Markov
chain, denoted $(\alb\cdot G, \mu)$, on the $G$-bimodule
$\alb\cdot G$ seen as a left $G$-space:
\[x_{n+1}\cdot g_{n+1}=h_{n+1}\cdot(x_n\cdot
g_n)=h_{n+1}(x_n)\cdot h_{n+1}|_{x_n}g_n.\]

If we start the Markov chain $(\alb\cdot G, \mu)$ at the point
$x=x\cdot 1\in\alb\cdot G$, we get a projection $\Pi_x:g_n\mapsto
g_n\cdot x=g_n(x)\cdot g_n|_x$ of the random walk on $G$ onto the
Markov chain on $\alb\cdot G$. If the action $(G, \alb)$ is
self-replicating (Definition~\ref{def:replicating}), then the map
$\Pi_x:G\arr G\cdot\alb:g\mapsto g\cdot x$ is onto.

The tuple $(\Pi_1(g), \ldots, \Pi_d(g))=(g(1)\cdot g|_1, \ldots,
g(d)\cdot g|_d)$ determines $g$ uniquely, since it means that the
image of $g$ in $\symm\wr_\alb G$ is $\pi(g|_1, \ldots, g|_d)$,
where $\pi$ is the permutation $i\mapsto g(i)$, $i\in\alb$.

The right action of $G$ on the bimodule $\alb\cdot G$ commutes
with the left action, hence the Markov chain $(\alb\cdot G, \mu)$
is invariant under the right action of $G$. The quotient of
$\alb\cdot G$ by the right $G$-action is naturally identified with
$\alb$ (the orbit of $x\cdot g$ is labeled by $x$). Consequently,
we get a Markov chain on $\alb$ equal to the quotient of the chain
$(\alb\cdot G, \mu)$ by the right $G$-action. It is easy to see
that this chain is given by the action of $G$ on $\alb$:
\[x\mapsto y,\qquad\text{with probability
$\sum_{g(x)=y}\mu(g)$.}\] This chain is irreducible when $G$ acts
transitively on the first level $\alb$ and $\mu$ is
non-degenerate.

Consider now the trace of the Markov chain $(\alb\cdot G, \mu)$ on
the right $G$-orbit $x\cdot G=\{x\cdot g\;:\;g\in G\}$ for a fixed
letter $x\in\alb$, which is a recurrent set (due to irreducibility
of the quotient chain), i.e., consider the Markov chain with the
first return transitions. We can identify the points of the orbit
$x\cdot G$ with the group $G$ using the map $x\cdot g\mapsto g$
and get in this way a Markov chain on $G$. It is not hard to check
that this Markov chain is the left random walk defined by the
measure $\mu_x=k_x(\mu)$, constructed above
(see~\cite{kaimanovich:munchhausen}).

\begin{example}
Let $B=\langle a, b\rangle$ be the Basilica group generated by the
states of the automaton shown on Figure~\ref{fig:basilica}. We
have $\phi(a)=\left(\begin{array}{cc} 1 & 0\\ 0 &
b\end{array}\right)$ and $\phi(b)=\left(\begin{array}{cc}0 & a\\
1 & 0\end{array}\right)$. Then $\phi(a^{-1})=\left(\begin{array}{cc} 1 & 0\\
0 & b^{-1}\end{array}\right)$ and
$\phi(b^{-1})=\left(\begin{array}{cc} 0 & 1\\ a^{-1} &
0\end{array}\right)$. Take
$\mu=p(a+a^{-1})+q(b+b^{-1})\in\ell_+^1(B)$, where $2(p+q)=1$.

The transition moves for the Markov chain on $\{0, 1\}\times B$
are
\begin{align*}
0\cdot g & \stackrel{a}{\arr} a(0)\cdot a|_0g=0\cdot g,\\
0\cdot g & \stackrel{a^{-1}}{\arr} a^{-1}(0)\cdot a^{-1}|_0g=0\cdot g,\\
1\cdot g & \stackrel{a}{\arr} a(1)\cdot a|_1g=1\cdot bg,\\
1\cdot g & \stackrel{a^{-1}}{\arr} a^{-1}(1)\cdot
a^{-1}|_1g=1\cdot b^{-1}g,
\end{align*}
and
\begin{align*}
0\cdot g & \stackrel{b}{\arr} b(0)\cdot b|_0g=1\cdot g,\\
0\cdot g & \stackrel{b^{-1}}{\arr} b^{-1}(0)=1\cdot a^{-1}g,\\
1\cdot g & \stackrel{b}{\arr} b(1)\cdot b|_1g=0\cdot ag,\\
1\cdot g & \stackrel{b^{-1}}{\arr} b^{-1}(1)\cdot
b^{-1}|_1g=0\cdot g.
\end{align*}
\end{example}

\begin{figure}
 \includegraphics{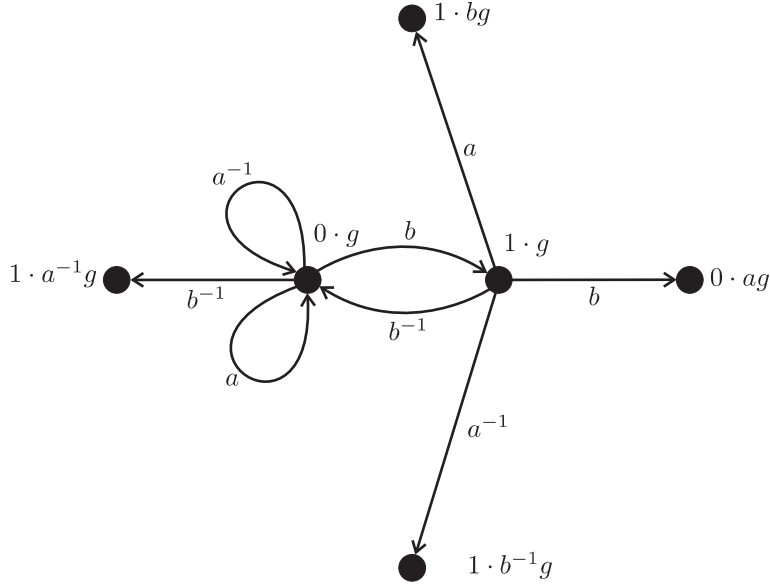}\\
 \caption{The random walk on $\{0, 1\}\cdot B$}\label{fig:graph}
\end{figure}

\begin{figure}
 \includegraphics{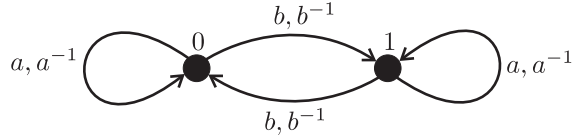}\\
 \caption{The induced random walk on $\alb$}\label{fig:smallgraph}
\end{figure}

See a graphical description of the random walk on $\{0, 1\}\cdot
B$ on Figure~\ref{fig:graph}. Figure~\ref{fig:smallgraph} shows
the graph of the induced random walk on $\{0, 1\}=\alb$ with
transition probabilities $1/2$.

The main characteristics of a random walk $(G, \mu)$ on groups
are:
\begin{itemize}
\item[(i)] the spectral radius of the random walk
\[r=\limsup_{n\to\infty}\sqrt[n]{\mu^{(n)}(e)};\]
\item[(ii)] the drift
\[\lambda=\lim_{n\to\infty}\frac{|g_n|}{n},\]
where $\{g_n\}_{n\ge 1}$ is a random trajectory and $|g_n|$ is the
length of an element;
\item[(iii)] entropy \[h=\lim_{n\to\infty}\frac 1n H(\mu^{(n)}),\]
where $H$ is the entropy of a probability distribution.
\end{itemize}

Vanishing of $\lambda$ implies vanishing of the entropy, which
implies amenability of $G$ (in case of a non-degenerate measure
$\mu$~\cite{kaimvershik}), while $r=1$ for a symmetric measure
$\mu$ is equivalent to amenability (see~\cite{kesten:amen}).

Kaimanovich's approach to amenability (called in his paper
``M\"unchhausen trick'') is based on the inequalities
(see~\cite{kaimanovich:munchhausen} Theorem~3.1)
\[h(G, \mu)\le h(G, \mu_i)\le |\alb|h(G, \mu)\]
which holds for the projections $\mu_i$, $i=1, \ldots, d$,
described above. His main observation based on these inequalities
and properties of the entropy is that if $\mu$ is non-degenerate
and self-affine, i.e., if there is $\alpha>0$ such that
\[\mu_i=(1-\alpha)\delta_e+\alpha\mu\]
for some $i$, then the inequalities imply $h(\mu)=0$ and hence
that $G$ is amenable (and moreover, has Liouville property).

Kaimanovich calls the measure $\mu$ satisfying the self-affinity
condition \emph{self-similar}. Self-similar measures are fixed
points of the maps
\[\alpha_i:\mu\mapsto\frac{k_i(\mu)-k_i(\mu)(e)}{1-k_i(\mu)(e)},\]
which are continuous maps of the simplex of probability measures
on $G$. The simplex is compact with respect to the weak topology.
For all $i$ there is an $\alpha_i$-invariant measure (indeed, the
measure $\delta_e$ concentrated in identity is such a measure) but
the problem is to find a non-degenerate fixed point if such exist.

The next few examples show when this is indeed the case. Here we
follow~\cite{kaimanovich:munchhausen}.

\begin{example}
Take $\mathfrak{G}=\langle a, b, c, d\rangle$ --- the group of
Example~\ref{sss:grigorchuk}. Consider a one-dimensional family of
measures
\[\mu=2\alpha m_1+4\beta m_2=\alpha a+\beta b+\beta c+\beta d+\alpha+\beta,\]
where $m_1=(1+a)/2$, $m_2=(1+b+c+d)/4$ and $2\alpha+4\beta=1$ and
$0\le \alpha, \beta\le 1$.

Then
\begin{align*}
\mu_1 &= \frac\alpha{1-\alpha}+4\beta
m_1+\frac{4\alpha\beta}{1-\alpha}m_2\\
\mu_2 &=
\frac\alpha{1-\alpha}+\frac{4\alpha\beta}{1-\alpha}m_1+4\beta m_2.
\end{align*}

In terms of the parameter $\alpha\in (0, 1/2)$ the corresponding
transformations $\phi_i$ take the form
\begin{align*}
\phi_1:\alpha &\mapsto
\frac{2\beta}{1-\frac\alpha{1-\alpha}}=\frac{1-\alpha}2\\
\phi_2:\alpha &\mapsto
\left.\frac{2\alpha\beta}{1-\alpha}\right/\left(1-\frac\alpha{1-\alpha}\right)=
\frac{2\alpha\beta}{1-2\alpha}=\frac\alpha 2.
\end{align*}
The only fixed point for $\phi_2$ corresponds to $\alpha=0$ and
$\mu$ is degenerate in this case, while for $\phi_1$ the value
$\alpha=1/3$ gives a fixed point corresponding to a non-degenerate
self-affine measure
\[\mu=\frac 23 m_1+\frac 13 m_2=\frac 5{12}+\frac 13 a+\frac
1{12}(b+c+d).\] Removing the atom at the identity we obtain a
self-affine measure concentrated on the generating set $\{a, b, c,
d\}$
\[\wt\mu=\frac 47+\frac 17(b+c+d)\] with self-similarity
coefficient $1/2$
\[\wt\mu_1=\frac 12+\frac 12\wt\mu\]
\end{example}

\begin{example}
Basilica group, $B=\langle a, b\rangle$ is defined by
\begin{alignat*}{2} a(0w) &= 0w, &\quad a(1w) &= 1b(w)\\
b(0w) &= 1w, &\quad b(1w) &= 0a(w),\end{alignat*} see
Example~\ref{sss:basilica}. Consider the one-parameter family of
measures
\[\mu=\frac{a+a^{-1}+rb+rb^{-1}}{2(r+1)},\quad r\ge 0.\]
Then
\begin{multline*}k_2(\mu)=\frac{b+b^{-1}}{2(r+1)}+\frac{r^2(2+a+a^{-1})}{4(1+r)^2}\left(1-\frac
1{r+1}\right)^{-1}\\
=\frac r{4(r+1)}(a+a^{-1})+\frac{b+b^{-1}}{2(r+1)}+\frac
r{2(r+1)}\\
=\frac{r+2}{2(r+1)}\left(\frac r{2(r+2)}(a+a^{-1})+\frac
1{r+2}(b+b^{-1})\right)+\frac r{2(r+1)}\\
=p(r)\phi_2(\mu)+q(r),
\end{multline*}
where $p(r)=\frac{r+2}{2(r+1)}$, $q(r)=\frac r{2(r+1)}$ and
\[\phi_2\left(\frac 1{2(r+1)}: \frac
r{2(r+1)}\right)=\left(\frac r{2(r+2)}:\frac 1{r+2}\right)\] is
the projectivization of the corresponding map, which represents
the map $\Lambda:\frac 1{2r}\mapsto \frac{r}2$ of $\R$ or
$z\mapsto\frac 1{4z}$ if $z=\frac 1{2r}$. This map has no fixed
points, but is periodic of period $2$; $\Lambda^2=id$. This fact
is used in~\cite{barthvirag} to prove amenability of $B$.

At the same time, as it is shown
in~\cite{kaimanovich:munchhausen}, the map $\phi_1$ has a fixed
point, which represents a non-degenerate measure and this gives
another way to prove amenability of $B$.
\end{example}

\begin{question}
Using the M\"unchhausen trick construct new interesting examples
of self-similar amenable but not elementary amenable groups.
\end{question}

\bibliographystyle{plain}
\bibliography{nekrash,mymath}

\end{document}